\documentclass[proc]{edpsmath}
\usepackage[T1]{fontenc}
\usepackage{algpseudocode, algorithm}
\usepackage{bbold}
\usepackage{graphicx}
\usepackage{url} 
\usepackage{mathtools}
\usepackage[colorlinks = true,
            linkcolor = blue,
            urlcolor  = blue,
            citecolor = blue,
            anchorcolor = blue]{hyperref}

\usepackage{subcaption}

\usepackage{xcolor}

\newcommand{\cB}{\ensuremath{\mathcal{B}}}

\newcommand{\cF}{\ensuremath{\mathcal{F}}}

\newcommand{\cM}{\ensuremath{\mathcal{M}}}

\newcommand{\cP}{\ensuremath{\mathcal{P}}}

\newcommand{\cU}{\ensuremath{\mathcal{U}}}

\newcommand{\cY}{\ensuremath{\mathcal{Y}}}
\newcommand{\cZ}{\ensuremath{\mathcal{Z}}}

\newcommand{\bN}{\ensuremath{\mathbb{N}}}

\newcommand{\bR}{\ensuremath{\mathbb{R}}}

\newcommand{\rA}{\ensuremath{\mathrm{A}}}

\newcommand{\rU}{\ensuremath{\mathrm{U}}}

\newcommand{\rW}{\ensuremath{\mathrm{W}}}

\newcommand{\xmin}{\ensuremath{x_{\min}}}
\newcommand{\xmax}{\ensuremath{x_{\max}}}

\renewcommand{\Pr}{\ensuremath{\cP_2(\Omega)}}
\newcommand{\bary}{\ensuremath{\mathrm{Bar}}}

\newcommand\icdfo[1]{\operatorname{icdf}_{#1}}

\DeclareMathOperator*{\argmin}{arg\,min}

\DeclareMathOperator{\pdf}{pdf}
\DeclareMathOperator{\cdf}{cdf}
\DeclareMathOperator{\icdf}{icdf}
\DeclareMathOperator{\iicdf}{iicdf}

\newcommand{\cond}{\; :\;}
\newcommand{\opt}{\textrm{opt}}

\begin{document}

\title{Wasserstein model reduction approach for parametrized flow problems in porous media}
\author{Beatrice Battisti}\address{Politecnico di Torino \& Université de Bordeaux} \&
\author{Tobias Blickhan}\address{Max-Planck-Institute for Plasma Physics \& Technische Universit\"at M\"unchen} \&
\author{Guillaume Enchery}\address{IFPEN} \&
\author{Virginie Ehrlacher}\address{Ecole Nationale des Ponts et Chaussées \& INRIA} \&
\author{Damiano Lombardi}\address{INRIA} \&
\author{Olga Mula}\address{Eindhoven University}
%
%
%
\begin{abstract} 
The aim of this work is to build a reduced-order model for parametrized porous media equations. The main challenge of this type of problems is that the Kolmogorov width of the solution manifold typically decays quite slowly and thus makes usual linear model-order reduction methods inappropriate. In this work, we investigate an adaptation of the methodology proposed in~\cite{ehrlacher2020nonlinear}, based on the use of Wasserstein barycenters~\cite{AC2011}, to the case of non-conservative problems. Numerical examples in one-dimensional test cases illustrate the advantages and limitations of this approach and suggest further research directions that we intend to explore in the future. 
\end{abstract}
\begin{resume}
Le but de ce travail est de construire un modèle réduit pour des problèmes d'écoulements en milieux poreux paramétrés. La difficulté principale de ce type de problèmes est que la distance de Kolmogorov de l'ensemble de solutions décroît lentement, rendant ainsi les méthodes de réduction de modèles linéaires usuelles inefficaces. Ici, nous proposons une adaptation de la méthodologie proposée dans~\cite{ehrlacher2020nonlinear}, utilisant des barycentres de Wasserstein~\cite{AC2011}, au cas de problèmes non conservatifs. Des tests numériques en dimension 1 permettent d'illustrer les avantages et les limitations de cette approche et d'identifier des pistes de recherche que nous souhaiterons aborder dans un futur travail. 
\end{resume}
\maketitle

\section*{Introduction}

Model-order reduction methods have proved to be extremely useful tools in order to accelerate parametric studies in various contexts. Given a Partial Differential Equation (PDE) involving parameters, the main goal of model-order reduction is to provide a fast approximation of the parameter-to-solution map. Classical strategies are linear in the sense that they rely on approximating the set of solutions with linear or affine spaces. This approach is the backbone of most existing methods among which stand the reduced basis method (see \cite{HRS2015,QMN2016}), the empirical interpolation method and 
its generalized version (G-EIM, see \cite{BMNP2004,GMNP2007,MM2013,MMT2016}), Principal Component Analysis (PCA, see \cite[Chapter 1]{BCOW2017}), polynomial-based methods like \cite{CDS2010,CDS2011} or low-rank methods \cite{KS2011}. The efficiency of such approaches is linked to the decay of the so-called Kolmogorov width of the solution manifold, which informs about the best possible approximation performance of the set of solutions when working with linear spaces. In cases where this width decays quickly with the dimension of the reduced space, significant gains in computational time can be achieved with respect to working with a direct, full-order computation of the solutions with classical discretization methods such as finite elements, or finite volumes. Families of PDEs where the Kolomogorov width presents a fast decay are elliptic or parabolic problems (see \cite{CD2015}).

However, in situations where the decay of the Kolmogorov width is slow with respect to the dimension of the reduced space, linear model-order reduction methods are inherently doomed to fail. This is the case, for instance, of parametrized flow problems in porous medium equations, which is the object of the present work. The acceleration of parametric studies of this type of problems requires the use of \itshape nonlinear \normalfont model-order reduction methods. This direction is a very active research topic, and numerous nonlinear strategies have been proposed (see \cite{carlberg2015adaptive, amsallem2012nonlinear,amsallem2016pebl, peherstorfer2014localized, lee2018model,gonzalez2018learning} to name a few). In this work, we consider an adaptation of a methodology based on the use of Wasserstein barycenters, which was originally proposed in~\cite{ehrlacher2020nonlinear} for the reduction of conservative flow problems (see also \cite{bonneel_wasserstein_2016}, \cite{schmitz_wasserstein_2018} for related work in computer graphics and computational anatomy). In our present context, we wish to investigate methods for the reduction of parametrized flow problems in porous media problems, which are not conservative. The fact that the total mass of each fluid phase in the system of interest may vary with respect to time and other parameter values yields a difficulty compared to the cases which were studied in~\cite{ehrlacher2020nonlinear}, where only conservative flows were considered. We propose here an approach in order to circumvent them. In this preliminary work, we only investigate one-dimensional problems, and leave the extension of the approach for higher-dimensional problems for future research. 

The paper is structured as follows: first, we present the parametrized two-phase flow problem in porous media which we focus on in this work in Section~\ref{sec:application}. Then, we recall some basic properties of Wasserstein spaces and Wasserstein barycenters in dimension 1 in Section~\ref{sec:nota} and we describe the proposed model-order reduction methodology in Section~\ref{sec:methodology}. Finally, numerical results illustrating the strengths and weaknesses of our approach are shown in Section~\ref{sec:numerical}.

\section{Problem Setting: Two-phase flow problem in porous media}
\label{sec:application}


\subsection{The physical model}
We consider a simplified model of two-phase flow through a porous rock. This type of model is often used in geosciences to describe the displacement of two non-miscible fluid phases. The wetting phase often corresponds to the water component whereas the second one depends on the application: air in hydrogeology, $\text{CO}_2$ for underground gas storage...

In the following, we denote  by $\Omega \subset \bR^d$, $d\in \bN^*$ the spatial domain occupied by the porous medium. The domain $\Omega$ has a total volume
$$V_t = V_v + V_s > 0$$
which is the sum of the volume $V_v$ occupied by the void, and the volume $V_s$ occupied by the solid material. We can thus characterize the medium by its porosity
$$
\phi\coloneqq\frac{V_v}{V_t} \in [0, 1]
$$
which is the ratio of the volume of void space $V_v$ over the total volume of the material $V_t$. {We assume here that the porosity of the medium is constant in time and characterized by a function $\phi: \Omega \to [0,1]$. } The volume $V_v$ of the void space is here assumed to be completely saturated by two incompressible and non-miscible fluid phases, called wetting ($\rm w$) and non-wetting (${\rm nw}$), with respective volumes $V_{\rm w}$ and $V_{{\rm nw}}$ which are such that $$V_v = V_{\rm w} + V_{{\rm nw}}.$$
Their volume fraction is given by the saturation
$$s_{\alpha}\coloneqq V_\alpha/V_v, \quad \alpha \in \{{\rm w}, {\rm nw}\}$$
and we have that
$$
s_{{\rm nw}}+s_{{\rm w}}=1.
$$
Due to  incoming/outcoming flows from outside the domain $\Omega$ (runoff waters from the surface, injecting/pumping wells...), both saturations $s_{nw}$ and $s_{w}$ depend on space and time. For a given time interval $[0, T]$, they are solutions to an equation stating the conservation of their phase volume,
\begin{equation}
    \phi(x) \ \partial_t s_{\alpha}(x,t) + \nabla \cdot  v_\alpha(x, t) = 0,
    \qquad \forall (x,t) \in \Omega \times (0, T),\;  \forall \alpha \in \{{\rm w}, {\rm nw}\}.
\label{eq:phase_conservation}
\end{equation}
In equation \eqref{eq:phase_conservation}, $v_\alpha$ denotes the generalized Darcy velocity. Neglecting the gravity and capillary effects, the velocity is defined for all $(x,t)\in \Omega\times (0, T)$ as
\begin{equation}
    v_\alpha(x, T) \coloneqq - \lambda_{\alpha}(x, t) \ {K}(x) \ \nabla p(x, t), \quad \forall \alpha \in \{{\rm w}, {\rm nw}\},
\label{eq:darcy}
\end{equation}
where $p$ is the fluid pressure. The pressure depends on space and time, and it is assumed to be the same for both phases. The coefficient {$K: \Omega \to \bR_+$} is the absolute permeability of the porous medium and reflects the ability of one fluid to flow through that medium. The phase mobility $\lambda_\alpha$ is also involved in the definition of the Darcy velocity. It is defined as
$$
\lambda_\alpha(x, t) \coloneqq \frac{k_{r_{\alpha}}(x, t)}{\mu_{\alpha}}
$$
where $\mu_{\alpha}\in \bR_+$ is the phase viscosity. $k_{r_{\alpha}}$ is the relative permeability which is an increasing function of $s_\alpha$ with values in $[0,1]$. In hydrogeology or in reservoir engineering applications, this law is often defined thanks to the Brooks-Corey model \cite{BROOKS}, which we simplify here to a power law
$$
k_{r_{\alpha}}(x, t)=s_{\alpha}^{\beta}(x, t),
$$
with $\beta>0$. The dependency of $k_{r_\alpha}$ on $s_\alpha$ induces a nonlinear coupling between equations \eqref{eq:phase_conservation} and equation \eqref{eq:darcy} for the unknowns are $s_{\rm w},\ s_{{\rm nw}}$ and $p$. In this work,  the system is closed with  boundary conditions
\begin{alignat}{3}
\label{eq:pressBC}
p(t, x) &= p_D(t, x),\quad &&\forall (x, t) \in \partial \Omega\times [0, T], \\
\label{eq:satBC}
s_{\rm w}(t, x) &= 1, \quad &&\forall (x, t) \in \partial \Omega\times [0, T] \;\text{ s.t. }\; v(t, x) \cdot n(x) < 0,
\end{alignat}
where $n$ is the outward unit normal of $\partial \Omega$, and initial conditions
\begin{equation}
s_{{\rm w}}(x,t=0) = s_{\rm w}^0,\quad \forall x\in \Omega. 
\label{eq:IC}
\end{equation}
In general, the problem is often reformulated in terms of pressure $p$ and the saturation $s_w$. By introducing the total Darcy velocity
\begin{equation}
\label{eq:totalVel}
v \coloneqq v_{\rm w} + v_{{\rm nw}} = -(\lambda_{\rm w} + \lambda_{\rm nw})k \nabla p,
\end{equation}
the previous system of equations is equivalent to
\begin{equation}
\label{eq:full-system}
    \begin{cases}
    &\textrm{div}(v) = 0 \\
    &\phi \ \partial_t s_{\rm w} + \nabla \cdot (f_{\rm w}(s_{\rm w}) v) = 0 \\
    &s_{\rm nw} =1 - s_{\rm w},
    \end{cases}
\end{equation}
where
$$
f_\alpha(s_{\alpha})\coloneqq \frac{\lambda_{\alpha}(s_{\alpha})}{\lambda_{\rm w}(s_{\rm w})+\lambda_{\rm nw}(s_{\rm nw})},\quad \alpha \in \{{\rm w}, {\rm nw}\}
$$
is the fractional flow.

Problem \eqref{eq:full-system} with conditions \eqref{eq:pressBC}-\eqref{eq:satBC}-\eqref{eq:IC} is solved using a finite-volume discretization in space and an IMPES (Implicit for Pressure, Explicit for Saturation) scheme in time. This latter scheme corresponds to an Euler implicit scheme for the pressure in the first equation of \eqref{eq:full-system} and to an explicit scheme for the saturation in the second equation of \eqref{eq:full-system} where the phase mobilities in \eqref{eq:totalVel} are evaluated explicitly too.
A similar discrete system, made of a two-point scheme for the pressure gradient and an upwind scheme for the saturation, has already been considered in \cite{boyavalEtAl2017} to study the reductibility of the first equation of \eqref{eq:full-system} with respect to some rock and fluid properties. The only difference between the present discretization and that study lies in the time discretation used for the saturation equation, which was implicit in this previous work and could allow the use of larger time steps. We refer to that work and to the references cited therein for more details on these discretizations. The IMPES scheme enables to decouple the resolution of \eqref{eq:full-system}, starting from the pressure, updating the total Darcy fluxes and finally computing the saturations.

\subsection{Parametric Problem and Manifold set of solutions} \label{sec:manifold}

Let $\cY \subset \bR^p$ be a compact set of parameter values for some $p\in \mathbb{N}^*$. For all $y\in\cY$, let $K^y: \Omega \to \mathbb{R}_+$, $\phi^y: \Omega \to [0,1]$, $\beta^y\in \mathbb{R}_+^*$ and $\mu_\alpha^y \in\mathbb{R}_+ $ for $\alpha\in\{ {\rm w}, {\rm nw}\}$. Moreover, let us denote by $s^y = s_{\rm w}^y$ the solution $s_{\rm w}$ to system \eqref{eq:full-system} with data $K = K^y$, $\phi = \phi^y$, $\beta = \beta_y$ and $\mu_\alpha = \mu_\alpha^y$.

One may need, for uncertainty quantification for instance, to evaluate solutions to \eqref{eq:full-system} for many instances of parameter values $y\in \cY$, and thus to reduce the computational time to evaluate the parameter-to-solution map because solving \eqref{eq:full-system} may be too computationally demanding.

More precisely, we wish to approximate the set of solutions
\begin{equation}
\label{eq:manifold}
\cM \coloneqq \left\{ s^y(\cdot, t) \cond (y,t)\in \cY \times [0,T]\right\}.
\end{equation}
Note that time is treated as parameter in $\cM$. Thus defining $\mathcal Z:= \cY \times [0,T]$ and denoting by $s^z$ the function $s^y(t,\cdot)$ for all $z=(y,t)\in \mathcal Z$, the solution set $\mathcal M$ can be rewritten as
\begin{equation}
\label{eq:manifold2}
\cM \coloneqq \left\{ s^z\cond z\in \mathcal Z\right\}.
\end{equation}


In the following, we develop a model reduction approach to approximate elements of the solution set $\cM$ based on the use of Wasserstein barycenters. 
 For the sake of simplicity, we will assume from now on that $d=1$ with one-dimensional problems, for which certain computations in this space are greatly facilitated and that $\Omega = (\xmin, \xmax) \subset \bR\text{, with } -\infty \leq \xmin<\xmax \leq \infty$. Section \ref{sec:nota} presents the essential notions of Wasserstein spaces and barycenters, and Section \ref{sec:methodology} presents our model reduction algorithm.


\section{Wasserstein barycenters in dimension 1}\label{sec:nota}

In this section, we recall some basic properties of the $L^2$-Wasserstein space in dimension 1 and introduce the notion of Wasserstein barycenters.

\subsection{Definition of the $L^2$-Wasserstein space in one dimension}\label{sec:not}

Let $\Pr$ denote the set of probability measures on $\Omega\subseteq \bR$ with finite second-order moments. For all $u\in \Pr$, we denote
\begin{equation}
\operatorname{cdf}_u : \left\{
\begin{array}{ccc}
 \Omega &\to & [0,1]  \\
 x & \mapsto & \operatorname{cdf}_u(x) := \int_{\xmin}^x \,du\\
\end{array}\right.
\end{equation}
its cumulative distribution function (cdf), and 
\begin{equation}
\operatorname{icdf}_u: 
\left\{
\begin{array}{ccc}
[0,1] &\to & \Omega  \\
p & \mapsto & \operatorname{cdf}_u^{-1}(p):=\inf\{ x\in \Omega, \; \operatorname{cdf}_u(x) \geq p\}\\
\end{array} \right.
\end{equation}
the generalized inverse of the cdf (icdf). The $L^2$-Wasserstein distance is defined by
$$
W_2(u,v) :=\mathop{\inf}_{\pi \in \Pi(u,v)} \left( \int_{\Omega \times \Omega} (x-y)^2 \,d\pi(x,y) \right)^{1/2},\quad \forall (u,v)\in \Pr \times \Pr,
$$
where $\Pi(u,v)$ is the set of probability measures on $\Omega \times \Omega$ with marginals $u$ and $v$. In the particular case of one dimensional marginal domains, it can be equivalently expressed using the inverse cumulative distribution functions  as
\begin{equation}\label{EqFlatReduction}
W_2(u,v) = \|\operatorname{icdf}_u - \operatorname{icdf}_v\|_{L^2([0,1])}.
\end{equation}
The space $\Pr$ endowed with the distance $W_2$ is a metric space, usually called $L^2$-Wasserstein space (see \cite{Villani2003} for more details). 

\subsection{Wasserstein barycenters}\label{subsec:wasserstein_barycenters}

Let $n\in \bN^*$ and let
$$
\Sigma_n:=\left\{ (\lambda_1,\cdots,\lambda_n) \in [0,1]^n, \quad \sum_{i=1}^n \lambda_i = 1\right\}
$$
be the probability simplex of dimension $n-1$. 

For any family of probability measures $\rU_n = (u_i)_{1\leq i\leq n} \in \Pr^n$ and barycentric weights 
$\Lambda_n = (\lambda_i)_{1\leq i\leq n} \in \Sigma_n$, there exists a unique minimizer to
\begin{equation}
\label{eq:bary}
\bary(\rU_n, \Lambda_n) = \argmin_{v\in \Pr} \sum_{i=1}^n \lambda_i W_2(v,u_i)^2,
\end{equation}
which is the barycenter with respect to $\rU_n$ and $\Lambda_n$. The barycenter can be easily characterized in terms of its inverse cumulative distribution function since \eqref{eq:bary} implies that
\begin{equation}
\label{eq:bary-icdf}
\icdfo{\bary(\rU_n, \Lambda_n)} = \mathop{\rm argmin}_{f\in L^2([0,1])} \sum_{i=1}^n \lambda_i \|\icdfo{u_i} -f\|_{L^2([0,1])}^2,
\end{equation}
which yields
\begin{equation}
\label{eq:bary-icdf-2}
\icdfo{\bary(\rU_n, \Lambda_n)} =\sum_{i=1}^n \lambda_i \icdfo{u_i}.
\end{equation}

For any given function $u\in \Pr$, one can define a notion of best barycenter $b(u, \rU_n)\in \Pr$ approximating $u$ among the family of barycenters
$$
\cB(\rU_n) = \{ \bary(\Lambda_n, \rU_n) \cond \Lambda_n \in \Sigma_n \}.
$$
The best barycenter is characterized in terms of its optimal weights
\begin{equation}\label{eq:minbary}
\Lambda_n^{\rm opt}  = \argmin_{\Lambda_n \in \Sigma_n}  W^2_2(u, \bary(\rU_n, \Lambda_n)),
\end{equation}
which are sometimes called the \textit{barycentric weights}. This allows us to define the optimal barycenter as
$$
b(u, \rU_n) = \bary(\rU_n, \Lambda_n^{\rm opt}).
$$
It can be easily characterized in terms of its inverse cumulative distribution function. This is due to the fact that for all $\Lambda_n\in \Sigma_n$ and all $u\in \Pr$,
\begin{align*}
 W^2_2(u, \bary(\rU_n, \Lambda_n)) &= \left\| \icdfo{u} - \sum_{i=1}^n \lambda_i \icdfo{u_i}\right\|^2_{L^2([0,1])}\\
\end{align*}
so the computation of the optimal weights $\Lambda_n^{\rm opt}$ in problem \eqref{eq:minbary} is a simple cone quadratic optimization problem. Note that at least one minimizer $\Lambda_n^\opt$ exists but uniqueness is not guaranteed.\\


\section{The Greedy Barycenter algorithm}
\label{sec:greedy-barycenter}
We detail in this section the greedy barycenter algorithm we use here so as to build a reduced-order model to efficiently compute approximations of elements of $\mathcal M$.  

\subsection{Offline phase: greedy algorithm}

For all $s\in \mathcal M$, we denote by $m(s):= \int_{\Omega} s$ and by $$
\mathcal U:= \left\{ u^z:=\frac{s^z}{m(s^z)}, \, z\in \mathcal Z\right\}.
$$

From now on, we make the assumption that $\mathcal U \subset \Pr$ and denote by $f^z:= \icdf_{u^z}\in L^2(0,1)$ for all $z\in \mathcal Z$.

\medskip

The Greedy Barycenter algorithm is an iterative procedure which aims at selecting, after $n$ iterations of the algorithm, a family of parameters $(z_1,\cdots,z_n)\subset \mathcal Z$, and hence a family $\rU_n = (u^{z_i})_{1\leq i\leq n}\subset \mathcal U$, such that any element $u\in \mathcal U$ can be approximated with good accuracy in the Wasserstein norm by a Wasserstein barycenter $\bary(\rU_n, \Lambda_n)$ for some $\Lambda_n \in \Sigma_n$.

    Let ${\mathcal Z}_{\rm train} \subseteq \mathcal Z$ be a finite subset of $\mathcal Z$ and let $\mathcal U_{\rm train} \coloneqq \{ u^z \cond z \in \mathcal Z_{\rm train} \}$ and
    \begin{equation}
        \cF_{\rm train} := \{ f^z= \icdf_{u^z} \in L^2(0,1): z \in \mathcal Z_{\rm train} \}.
    \end{equation}
    In practice, the strategy consists in finding a good collection of \textit{atoms} $\rA_n = \left ( a_i := f^{z_i}\right )_{1\leq i \leq n}$ such that for all $z\in \mathcal Z$, $f^z = \icdf_{u^z}$ is approximated with good accuracy in the $L^2(0,1)$ norm by a convex combination of elements of $\rA_n$, following similar guidelines as in~\cite{ehrlacher2020nonlinear}. \\
    
    \medskip
    
    \bfseries Greedy Barycenter Algorithm: \normalfont
    
    \begin{itemize}
     \item \bfseries Initialisation: \normalfont find $(z_1,z_2)\in \mathcal Z_{\rm train}$ solution to 
     $$
     (z_1,z_2) \in {\mathop{\rm argmax}}_{z,z'\in \mathcal Z_{\rm train}} \left\| f^z - f^{z'}\right\|_{L^2(0,1)}^2.
     $$
     Set $\rU_2:=(u^{z_1}, u^{z_2})$ and $\rA_2:=(f^{z_1}, f^{z_2})$.
 \item \bfseries Iteration $n\geq 3$: \normalfont fin $z_n \in \mathcal Z_{\rm train}$ such that
 \begin{equation}\label{eq:atoms}
 z_n \in \mathop{\rm argmax}_{z\in \mathcal Z_{\rm train}} \mathop{\min}_{(\lambda_i)_{1\leq i \leq n-1} \in \Sigma_{n-1}} \left\| f^z - \sum_{i=1}^{n-1} \lambda_i f^{z_i}\right\|_{L^2(0,1)}^2.
 \end{equation}
   Set $\rU_n:=(u^{z_1}, u^{z_2}, \cdots, u^{z_n})$ and $\rA_n:=(f^{z_1}, f^{z_2}, \cdots, f^{z_n})$.  
     \end{itemize}
  
The obtained set of atoms $\rA_n = (f^{z_1}, \cdots, f^{z_n}) = (a_1, \dots, a_n)$ is the output of the greedy algorithm.

As stated in Section \ref{subsec:wasserstein_barycenters}, solving \eqref{eq:atoms} amounts to solving $|{\mathcal Z}_{\rm train}|$ quadratic optimization problems on unit simplices in $\mathbb R^n$. These problems are solved using the COSMO package \cite{garstka_cosmo_2020}. While the number of optimizations problems to be solved in each greedy iteration is large, they are independent of each other and could be easily parallelized. When solving them consecutively, we can use the same solver object for all of them, which saves significant computation time. More details on the optimization problem will be discussed in Section \ref{sec:comp_cost}.

	The greedy algorithm terminates with one of the following three termination criteria:

	\begin{enumerate}
	\item when an absolute tolerance $\varepsilon_{abs}\geq 0$ is reached, i.e.
	\begin{equation}
		\Delta^2_n:= \max_{ z \in \mathcal Z_{\rm train} } \min_{\mathbf (\lambda_i)_{1\leq i  \leq n} \in \Sigma_n} \left \Vert f^z - \sum_{i=1}^n \lambda_i f^{z_i} \right \Vert^2_{L^2(0,1)} < \varepsilon_{abs}^2, 
	\end{equation}
	\item when the addition of a supplementary atom to the dictionary does not reduce the error beyond a given relative tolerance $\varepsilon_{rel}$, i.e.
	\begin{equation}
		\Delta_n - \Delta_{n+1} < \varepsilon_{rel} \, \Delta_n, 
	\end{equation}
	\item when a specified maximum number $n_{\rm max}$ of atoms is reached.
	\end{enumerate}
	
	\subsection{Online phase: interpolation}
	
	Let us assume here that we know the values of $m(s^z)$ for all $z\in \mathcal Z_{\rm train}$. We then denote by $\mathcal I_m: \mathcal Z \to \mathbb{R}_+$ an interpolation function such that 
	$$
	\forall z\in \mathcal Z_{\rm train},\quad \mathcal I_m(z) = m(s^z).
	$$
	Furthermore, let us assume that for all $z\in \mathcal Z_{\rm train}$, we know the values $\Lambda^{\rm opt}(z) = (\lambda_1^{\rm opt}(z), \cdots, \lambda_n^{\rm opt}(z)) \in \Sigma_n$ solution to the optimization problem
	$$
	\Lambda_n^{\rm opt}(z) \in \mathop{\rm argmin}_{(\lambda_i)_{1\leq i \leq n}\in \Sigma_n}\left\| f^z - \sum_{i=1}^n \lambda_i a_i \right\|_{L^2(0,1)}^2.
	$$
	We then denote, for all $1\leq i \leq n$, by $\mathcal I_i: \mathcal Z \to \Sigma_n$ an interpolation function such that
	$$
	\forall z \in \mathcal Z_{\rm train},\quad \mathcal I_i(z) = \lambda_i^{\rm opt}(z).
	$$
	

The online phase of the reduced-order modeling method then works as follows:

\medskip

\bfseries Online phase: \normalfont
    \begin{enumerate}
        \item Given a new parameter value $z^*\in \mathcal Z$, we compute for all $1\leq i \leq n$, $\widetilde{\lambda}_i^* = \mathcal I_i (z^*)$ and $\widetilde{m}^* = \mathcal I_m (z^*)$. 
        \item We compute $\Lambda^*$ as the projection on $\Sigma_n$ of  $\tilde \Lambda^*:= (\tilde \lambda_1^*, \cdots, \tilde \lambda_n^*)$ using Algorithm 1 from \cite{wang_projection_2013} and define $m^* = [\tilde m^*]_+$.
        \item Compute $f^*:=\sum_{i=1}^n \lambda^*_i a_i$.
        \item By construction, $f^*$ is then equal to $\icdf_{u^*}$ for some unique $u^* \in \Pr$. Compute $u^*$. 
        \item The reduced-order modeling approximation $s_n^{{\rm gB},z^*}$ of $s^{z^*}$ is then defined as $
        s_n^{{\rm gB},z^*} = m^* u^*$.
    \end{enumerate}

\section{Computational aspects}\label{sec:methodology}

The aim of this section is to present the methodology used in order to build a reduced-order model for the parametrized porous media flow problem presented in the preceding section. We also discuss some choices we make in our implementation. The code used to generate the examples shown here is available at
\begin{center}
\href{https://github.com/ToBlick/MorporJ}{https://github.com/ToBlick/MorporJ}
\end{center}
The link also contains certain interactive examples to help gain understanding on the behavior of weighted Wasserstein barycenters as one varies the weights.

\subsection{Computing cumulative distribution functions and their inverse} 

	In practice, a snapshot $s\in \mathcal M$ (i.e. the wetting saturation $s_w$) is provided, for a given value of the set of parameters, as a piece-wise constant function on a uniform grid of size $N = 1000$. We denote this discrete representation as $\hat{\mathtt s} \in \bR^N$ and use typewriter font for all discrete distributions, cumulative distribution functions, and their inverse. \\
	
	In our one-dimensional numerical examples (with $x_{\rm min} =0$ and $x_{\rm max}=1$ for the sake of simplicity), we always consider that the saturation $s$ has to satify the boundary condition $s(x_{\rm min}) = 1$. This is the reason why, from $\hat{\mathtt s}$, we compute another vector $\mathtt{s}:= \left( \mathtt{s}_i \right)_{1\leq i \leq N+2}\in \mathbb{R}^{N+2}$ such that $\mathtt{s}_1 = 0.0$, $\mathtt{s}_2 = 1.0$ and $\mathtt{s}_{i+2} = \hat{\mathtt{s}}_{i}$ for all $1\leq i \leq N$. The value of $\mathtt{s}_2$ is chosen so that to enforce the boundary condition mentioned above, and the value of $\mathtt{s}_1$ is chosen so that the algorithm we use to compute a discrete approximation of the icdf of $u=\frac{s}{m(s)}$ is more numerically stable. \\

	Let us describe this procedure here. First, from $\mathtt{s}$, an approximation of $u$ is given by the vector $\mathtt{u} = (\mathtt{u}_i)_{1\leq i \leq N+2}\in \mathbb{R}^{N+2}$ defined by
	$$
	\forall 1\leq i \leq N+2, \; \mathtt{u}_i  = \frac{\mathtt{s}_i}{\sum_{j=1}^{N+2} \mathtt{s}_j}.
	$$
	Second, an approximation of $\cdf{u}$ is given by the vector $\mathtt{cdf} = (\mathtt{cdf}_i)_{1\leq i \leq N+2}\in \mathbb{R}^{N+2}$ defined by
	$$
	\forall 1\leq i \leq N+2, \quad \mathtt{cdf}_i = \sum_{j=1}^i \mathtt{u}_j. 
	$$
	Note that the definition of $\mathtt{s}$ guarantees that $\mathtt{cdf}_1 = 0$ and $\mathtt{cdf}_{N+2} = 1$.

    The discrete approximation of $\icdf_{u}$ will be given as a piecewise constant function on the interval $[0,1]$. Let $\mathtt{x}_1 = x_{\rm min} < \mathtt{x}_2 < \cdots < \mathtt{x}_{N+2} = x_{\rm max}$ be a uniform discretization grid of $(x_{\rm min}, x_{\rm max})$. Let also $M\in \mathbb{N}^*$ (in practice $M$ is chosen to be equal to $N+2$ in the numerical experiments) and $\mathtt{p}_1 = 0 < \mathtt{p}_2 < \cdots < \mathtt{p}_M =1$ be a uniform discretization grid of the interval $(0,1)$. The discrete approximation of $\icdf_{u}$ is then given as a vector $\mathtt{icdf} = (\mathtt{icdf}_j)_{1\leq j \leq M} \in \mathbb{R}^M$ the components of which are defined as follows:  for all $1\leq j \leq M$, let $i\geq 2$ be the smallest integer such that $\mathtt{cdf}_i \geq \mathtt{p}_j$. Then, 
    
    \begin{equation}
        \mathtt{ icdf}_j = \mathtt{x}_{i-1} + (\mathtt{x}_i - \mathtt{x}_{i-1}) (\mathtt{p}_j - \mathtt{cdf}_{i-1})/(\mathtt{cdf}_i - \mathtt{cdf}_{i-1}).
    \end{equation}
    
    Note that this can be done in a single pass over \texttt{j}, since the cdf is monotonically increasing. \\
	
	If we apply the same procedure to the resulting icdf, with the roles of \texttt{x} and \texttt{p} swapped, we obtain a numerical approximation of the inverse of the inverse cumulative distribution function, which is given as a vector $\mathtt{iicdf}\in \mathbb{R}^{N+2}$. Up to the error introduced by the finite grid size and the linear interpolation, $\mathtt{iicdf} \approx \mathtt{cdf}$. On average, this error lies between $10^{-4}$ and $10^{-3}$ in the discrete $L^2$ norm for the data from our numerical experiments. \\

	\medskip
	
	As mentioned at the end of Section \ref{subsec:wasserstein_barycenters}, the solution to the minimization problem yields the icdf of a normalized snapshot. We invert it to obtain the iicdf and perform simple first order backwards finite differences to compute the corresponding pdf. This guarantees that the reconstructed pdf integrates to one.

	Clearly, both the inversion and differentiation could be carried out using higher-order methods such as a spline interpolant of the discrete pdf. We choose the simple linear interpolants as they prove to be robust when applied to the snapshots we are considering in this work, which feature jump discontinuities and, in some cases, resemble even discrete Dirac distributions. This, of course, limits the numerical accuracy of our method, as will be discussed later.

\subsection{Sources of Error}
    Given a target $s^* = s^{z^*}$ for some $z^*\in \mathcal Z$ and the reduced-order model approximation obtained $s_n^{{\rm gB} \, *} = s_n^{{\rm gB},z^*}$, we can identify three sources of error.
    \begin{enumerate}
    \item 
    Firstly, $s^*$ might not be close to any $s_n^{{\rm gB}}$ that can be obtained by means of inverting and differentiating $\sum_{i=1}^n \mathcal I_i (z^*) f^{z_i}$ using the given atoms and interpolations. This error describes how rich the range of the collection of atoms is. The decay of this error as $n$ increases is of prime interest to evaluate the quality of the proposed method.
    \item Secondly, the optimization may fail to find the optimal barycentric weight $\Lambda_n^\opt$. The properties of the system that has to be solved is addressed in Section \ref{sec:comp_cost}.
    \item Thirdly, there is a numerical error introduced simply by following the $\pdf \mapsto \cdf \mapsto \icdf \mapsto \iicdf \mapsto \pdf$ pipeline as outlined in the beginning of this section. Naturally, this sets a lower bound on the error we can obtain, even with arbitrary large $n$. For all our numerical experiments, however, this error is at least one order of magnitude smaller than the average errors we observe in the numerical experiments, so it is not the limiting factor in these cases.
    \end{enumerate}
    Note that the proposed method operates entirely in the $W_2$ space and on the icdfs of the saturations. This error in $W_2$ need not reflect the errors in $L^1$ between the reconstructed saturations themselves, which  we chose as a reference. \\
    For example, consider the distributions $s_1 = 10 \times \mathbb{1}_{[0.0,0.1]}, s_2 = 2 \times \mathbb{1}_{[0.0,0.5]}$ and $s_3 = \mathbb{1}_{[0.0,1.0]}$. While $\Vert s_1 - s_2 \Vert_{L_1} = 1.6 > 1.0 = \Vert s_2 - s_3 \Vert_{L_1}$, we find $W_2(s_1,s_2) \approx 0.23 < 0.29 \approx W_2(s_2,s_3)$. Note that the distributions in this example are shaped similarly to the ones in our numerical experiments. As a result, reducing the training error in $W_2$ through the addition of an atom to the dictionary through the greedy algorithm may not reflect in the $L_1$ norm. \\
    The properties of the interpolation function are not adding to the error in our experiments, as we are only looking at the training error. This decision is made to keep the presentation brief as the generalization error depends largely on the selection of the training and test data sets, which we do not address.
    
\subsection{Short comments on  computational cost}\label{sec:comp_cost}

    A thorough study of the run-times for the proposed method is beyond the scope of this article. The numerical inversion of a discrete cdf or icdf depends only linearly on $N$, which allows the computation of Wasserstein barycenters for given weights in real time, as can be seen in the \texttt{demo} notebook accompanying the code provided along with this work.
    
    Selecting atoms with the greedy algorithm involves solving quadratic systems of the form $\min \Vert \mathtt A \mathtt \Lambda - \mathtt f \Vert^2$ with constraints $\mathtt{\Lambda \geq 0}$ element-wise and $\mathtt{\sum_i \lambda_i = 1}$. $\mathtt A$ is an $N \times n$ matrix with the discrete atoms as its columns, thus $n$ increases with every iteration of the algorithm. $\mathtt \Lambda$ is the vector of barycentric weights to be optimized and $\mathtt f$ is the icdf of the saturation we aim to approximate. As stated, this system has to be solved $\vert \cM_{\rm train} \vert$ times in every iteration and we observe that $\mathtt{A^T A}$ becomes very ill-conditioned with growing iteration number, see Figure \ref{fig:conditioning} below.
    We cannot orthogonalize the columns of $\mathtt{A}$, since the resulting atoms would no longer be discrete icdfs (for example, they will most likely no longer be strictly positive). This point is discussed further in Section \ref{sec:conditioning_convexity}.
     
    Since $\mathtt A$ does not change within greedy iterations, significant computational time can be saved by assembling the system once and then modifying $\mathtt f$ in-place. Furthermore, we initialize $\mathtt \Lambda$ inversely proportional to the $W_2$ distance between each atom in the dictionary and the distribution $\mathtt f$ that is to be approximated - that is, if $\mathtt f$ is very close to the i-th atom, the corresponding $\mathtt{\lambda[i]}$ will be initialized close to one.

\section{Numerical Examples}\label{sec:numerical}

    In this section, we compare the performance of the method for two test-cases to a POD representation. In the latter case, we project every element of $\cM_{\rm train}$ to the POD basis and back. We choose the relative $L^1$ error to compare the methods. To be precise, let $\mathtt \Psi$ be the $N \times n$ matrix that has as columns the first $n$ POD modes. The POD reconstruction of a solution $s$ is then $\mathtt{s^{POD} := \Psi_{POD}^T \Psi_{POD} s}$ and the corresponding error $\Vert s - s^{POD} \Vert_{L_1} / \Vert s \Vert_{L_1}$. The reconstruction error for the barycentric method is given by $\Vert s - \bar s^{gB} \Vert_{L_1} / \Vert s \Vert_{L_1}$, where $\bar s^{gB}$ is obtained by applying the method as it was described in Section \ref{sec:nota}. \\
    In the following examples, $\Omega=(0,1)$ , corresponding to a physical domain of length $1 \, \mathrm{km}$ and the grid is composed of $N+2=1002$ cells.
    The time interval, which covers several years, is discretized thanks to time steps satisfying the CFL condition given in \cite{EncheryMassonWolfEymard2002}. The boundary and initial conditions are: $p_D(0, t)=4.137 \times 10^7 \, \mathrm{Pa}$, $p_D(L, t)=2.758 \times 10^7 \, \mathrm{Pa}$, $s(0, t)=1$ for all $t \in [0, T]$, $s_w(x,0)=0$ for all $x \in (0, L)$.

\subsection{Example 1: Homogeneous medium with a varying exponent of the relative permeability law and various viscosity ratios} 

\begin{figure}[htbp]
        \centering
        \includegraphics[width=0.5\textwidth]{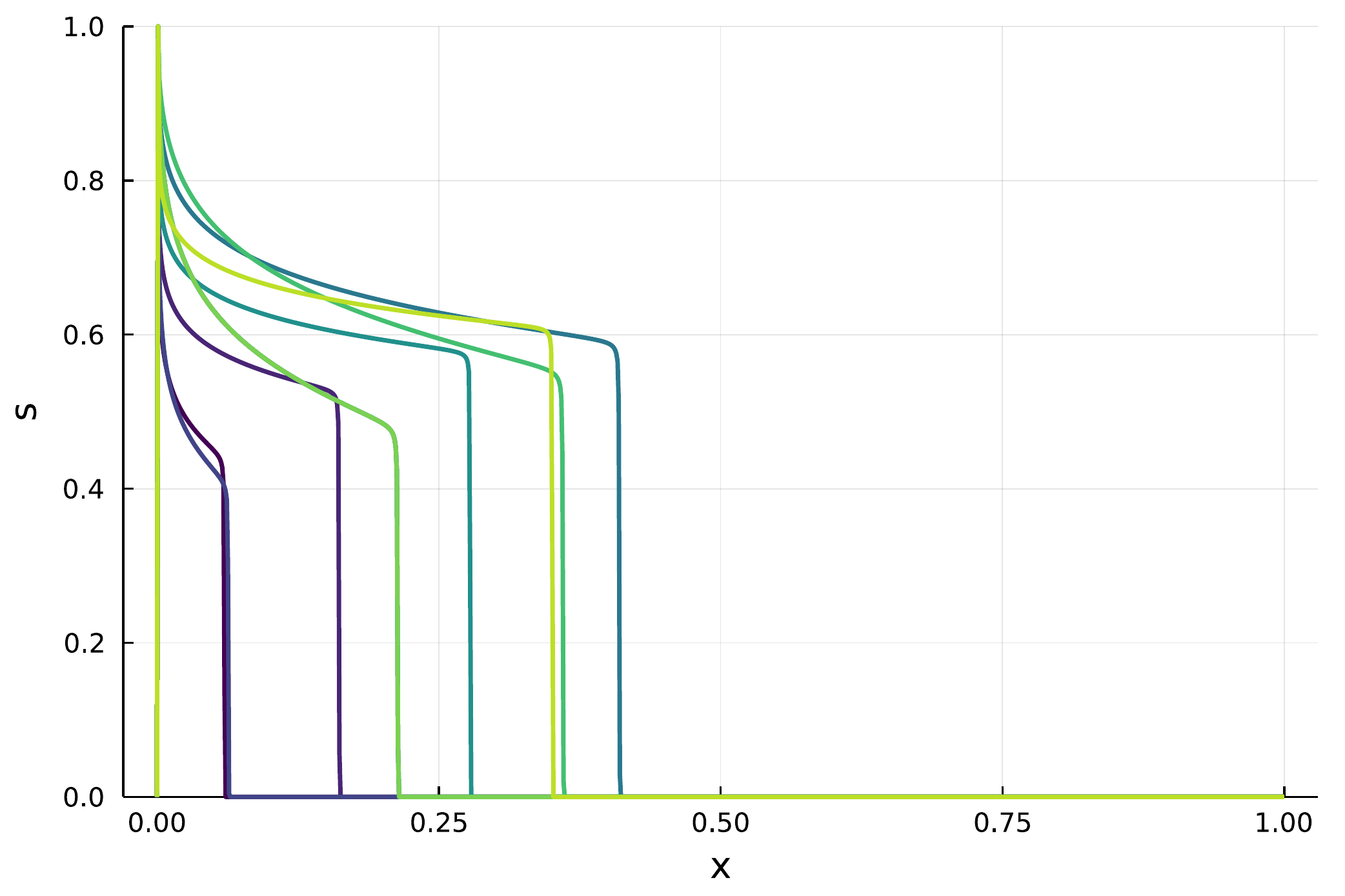}
        \caption{Randomly selected snapshots at final time from example 1.}
        \label{fig:ex1_snaps}
    \end{figure}

    For the first example we consider the case of a homogeneous medium with porosity $\phi=0.1$ and permeability $k= 10^{-13} \, \mathrm{m^2}$. The wetting viscosity is kept constant at a value of $\mu_{nw}=0.003 \, \mathrm{Pa\,s}$, but the non-wetting viscosity is changed in such a way that $\mu= {\mu_{nw}} / {\mu_{w}}$ ranges from 1 to 25. The sampling is chosen analogous to a geometric progression, so that $\mu \in \{ 1,2,3,6,12,25 \}$. The exponent $\beta$ of the relative permeability varies uniformly from 2 to 6, that is $\beta \in \{2,3,4,5,6\}$. The evolution of the water saturation is calculated over 5 years, saving {snapshots} at a fixed time-step, $t \in \{0.2, 0.4, \dots, 5.0 \}$ years, leading to $|  \cZ_{\rm train} | = 750$ and $z = (t,\mu,\beta)$. Some snapshots of the training set at the final time are shown in Figure \ref{fig:ex1_snaps}.
    
\begin{figure}[htbp]
    \centering
    \begin{subfigure}{.47\textwidth}
        \centering
        \includegraphics[width=0.92\textwidth]{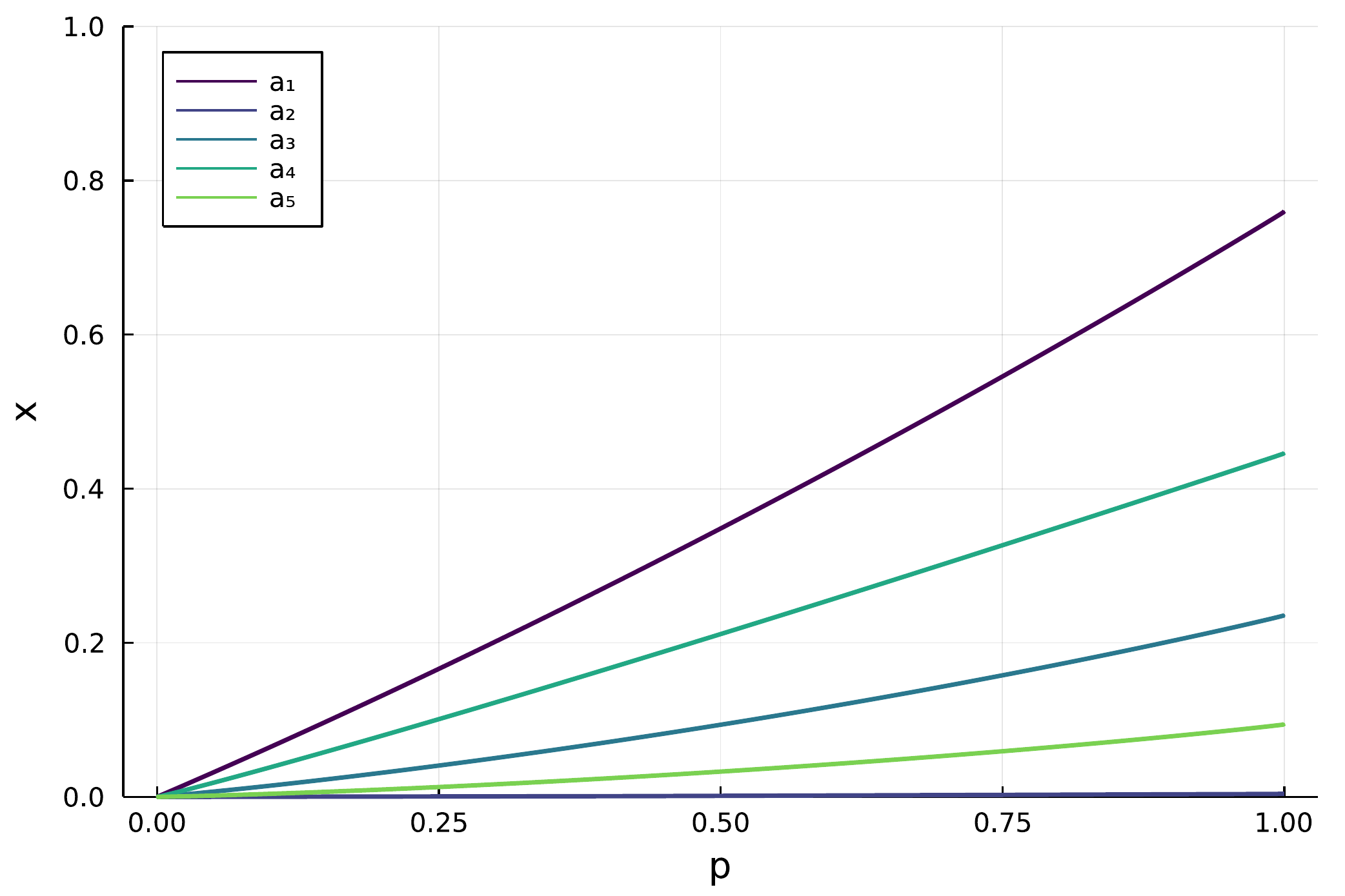}
        \caption{Dictionary atoms...}
    \end{subfigure}%
    \hspace{0.04\textwidth}
    \begin{subfigure}{.47\textwidth}
    \centering
        \includegraphics[width=0.92\textwidth]{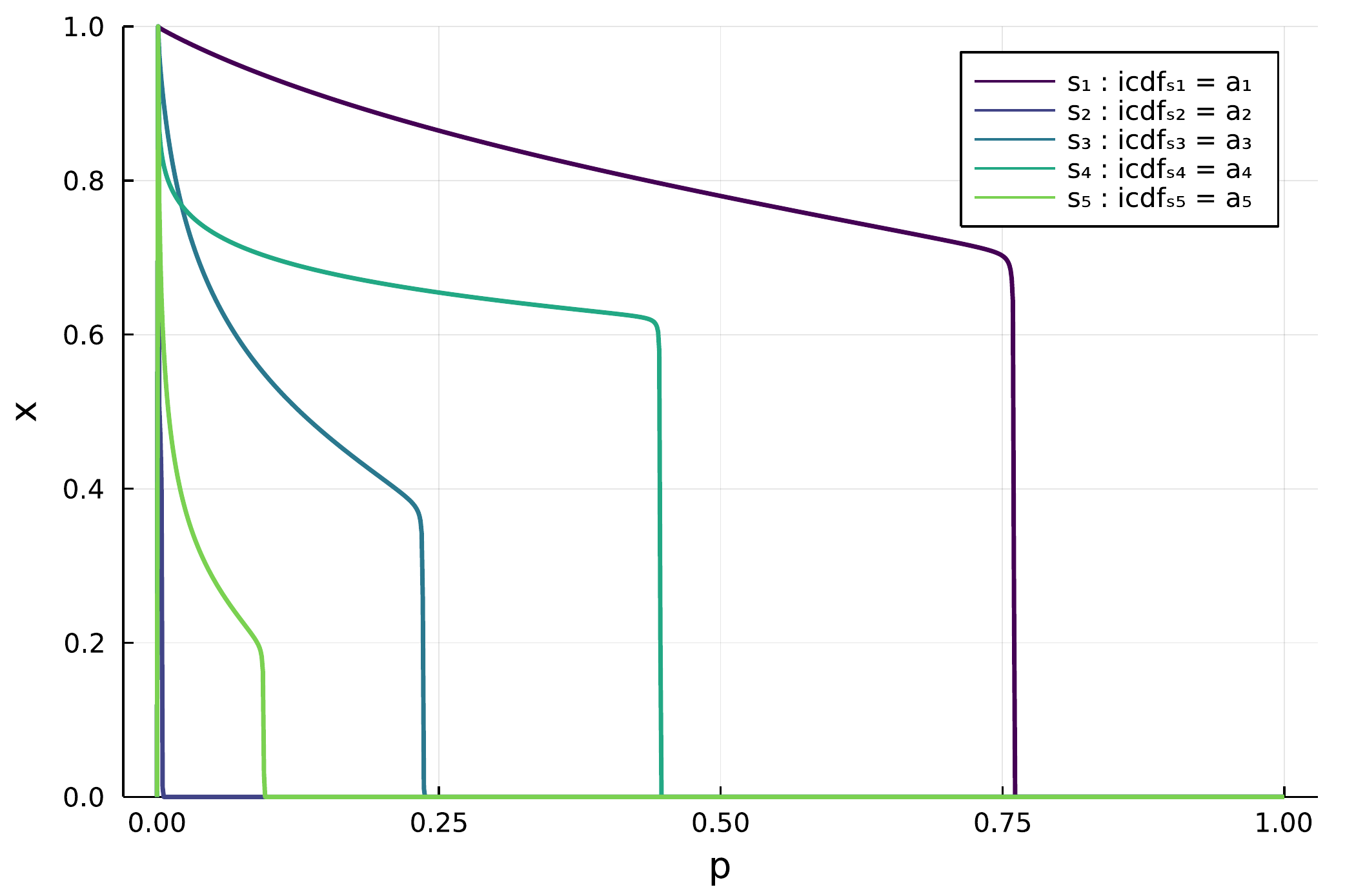}
        \caption{...and corresponding densities.}
    \end{subfigure}
    \caption{First five atoms chosen by the greedy algorithm for example 1.  }
    \label{fig:ex1_atoms}
    \end{figure}

\begin{figure}[htbp]
    \centering
    \begin{subfigure}{.47\textwidth}
    \centering
    \includegraphics[width=0.92\textwidth]{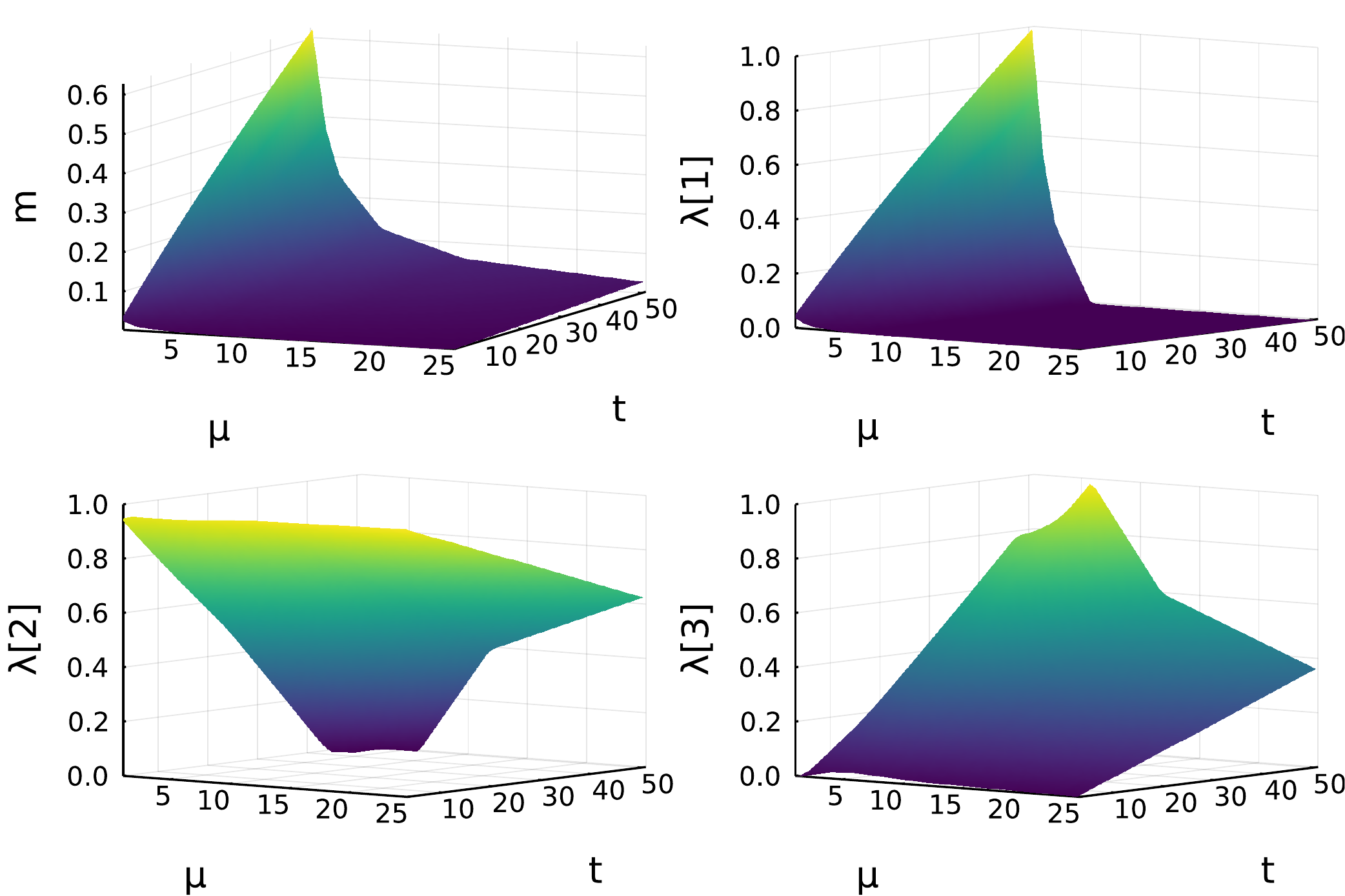}
    \caption{$(\mu,t) \mapsto \lambda_{1 \leq i \leq 3}(t, \mu, \beta = 2)$ and $(\mu,t) \mapsto m(t, \mu, \beta = 2)$.}
    \label{fig:ex1_interpolates1}
    \end{subfigure}%
    \hspace{0.04\textwidth}
    \begin{subfigure}{.47\textwidth}
    \centering
    \includegraphics[width=0.92\textwidth]{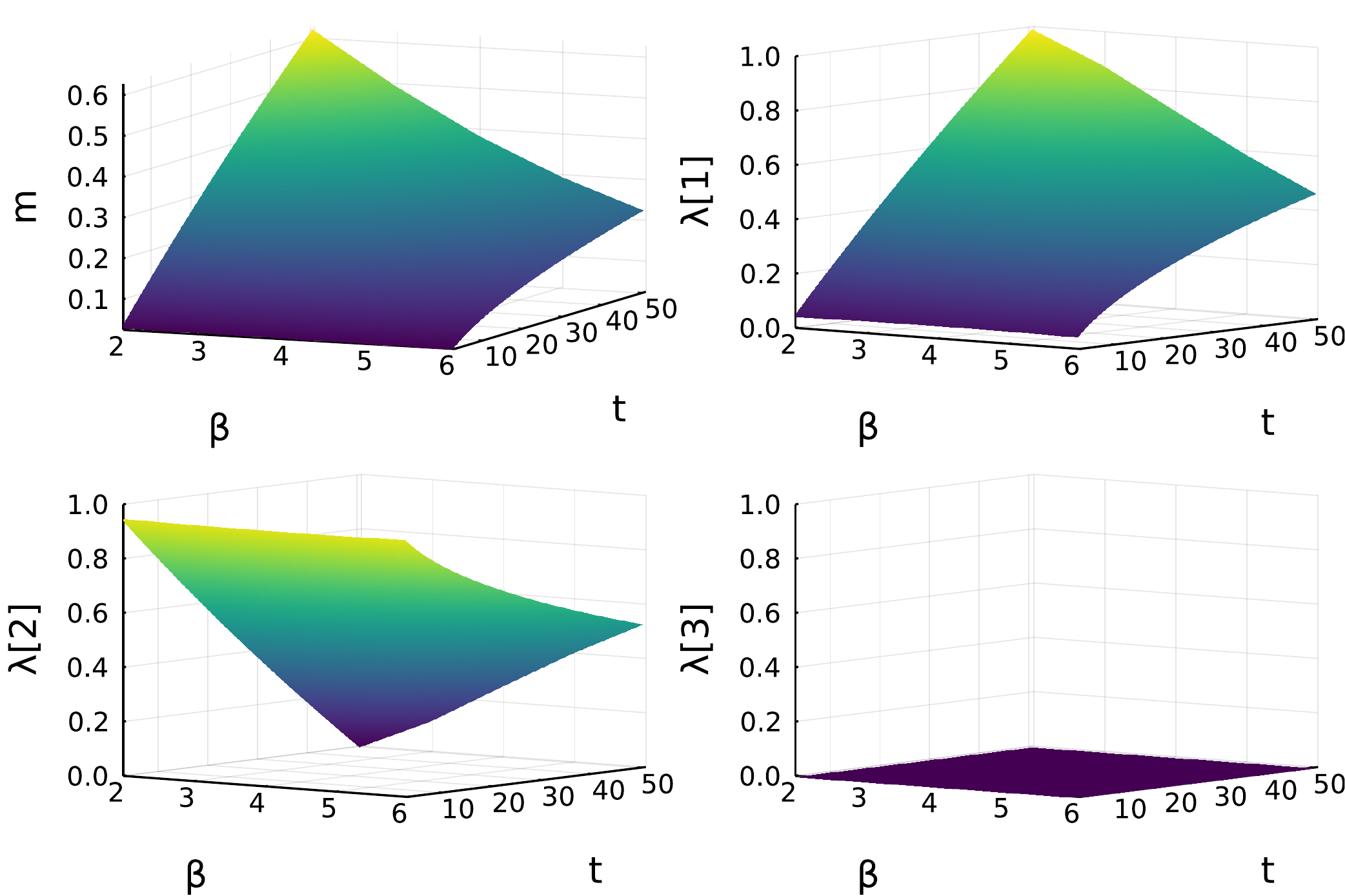}
    \caption{$(\beta,t) \mapsto \lambda_{1 \leq i \leq 3}(t, \mu=1, \beta)$ and $(\beta,t) \mapsto m(t, \mu=1, \beta)$.}
    \label{fig:ex1_interpolates2}
    \end{subfigure}
    \caption{Interpolations from parameter space to barycentric weights and mass in Example~1.}
    \label{fig:ex1_interpolates}
    \end{figure}

    The methodology described in Sections \ref{sec:greedy-barycenter}-\ref{sec:methodology} is then applied. A few selected atoms selected by the greedy algorithm are shown in Figure \ref{fig:ex1_atoms}. Figure \ref{fig:ex1_interpolates} shows cross-sections of the interpolating functions from the parameter space for the case $n=3$. Note how these are simple functions showing a smooth mapping from parameter space to the reduced coordinates. The average reconstruction error, measured in the $L_1$ norm, is already below 5\% for $n=3$ (see Table \ref{tab:ex1_n_of_epsilon}).
    
\begin{figure}[htbp]
        \centering
        \includegraphics[width=0.5\textwidth]{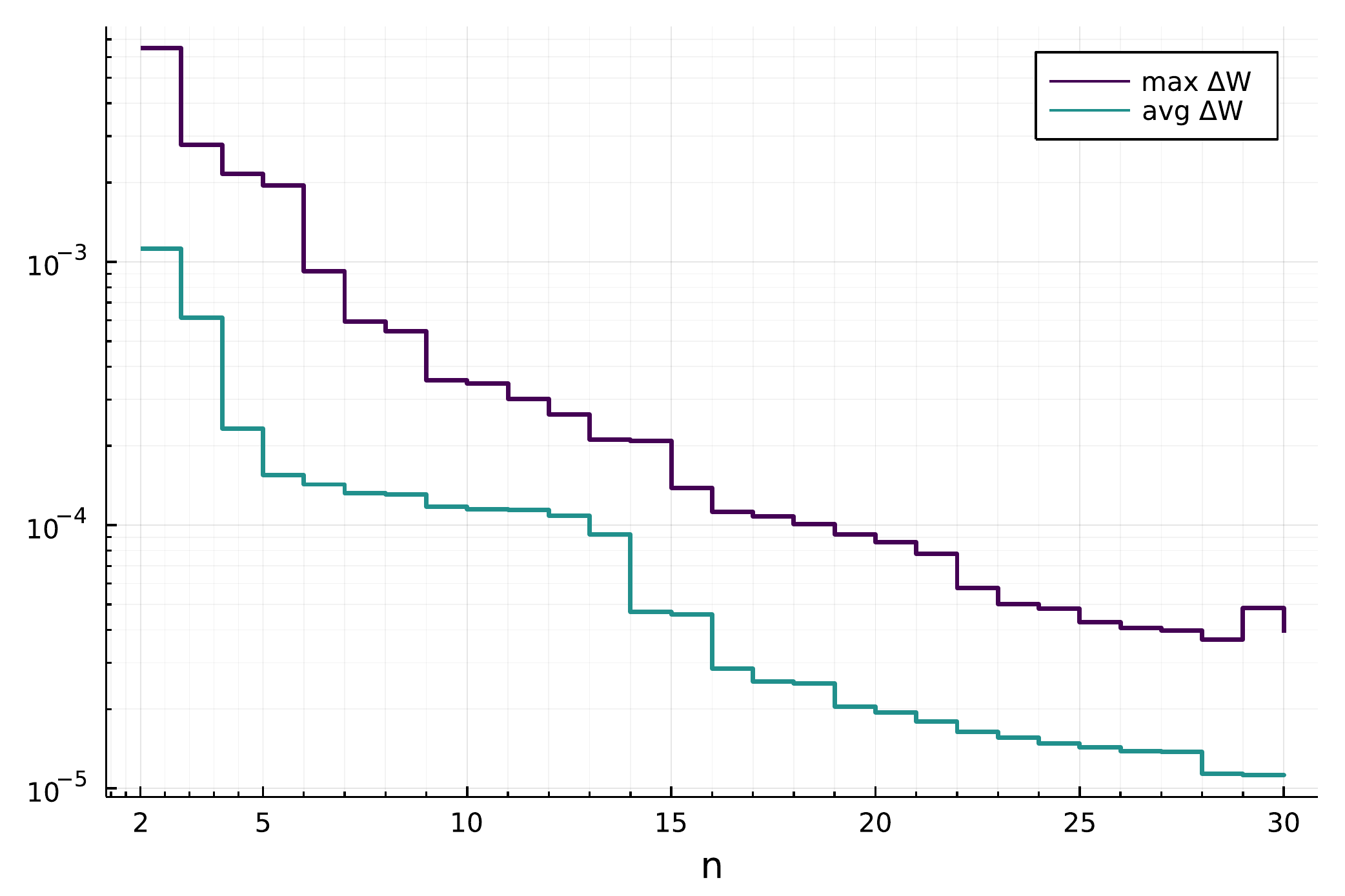}
        \caption{Decay of the $W_2$ error through the iterations of the greedy algorithm for example 1.}
        \label{fig:ex1_deltaW}
    \end{figure}

    The decay of the maximum and average Wasserstein error on $\cU_{\rm train}$ as more atoms are added is shown in Figure \ref{fig:ex1_deltaW}. Note that there is an increase in maximum error at the addition of the 29th atom. This can be explained by the iterative solver not converging before reaching its maximum iteration number while determining the optimal weights for one of the training snapshots. While the addition of atoms should only decrease the reconstruction error, we observe in practice that many atoms may make it harder to solve for the optimal barycentric weights (c.f. Section \ref{sec:conditioning_convexity}).
    
\begin{figure}[htbp]
        \centering
        \includegraphics[width=0.5\textwidth]{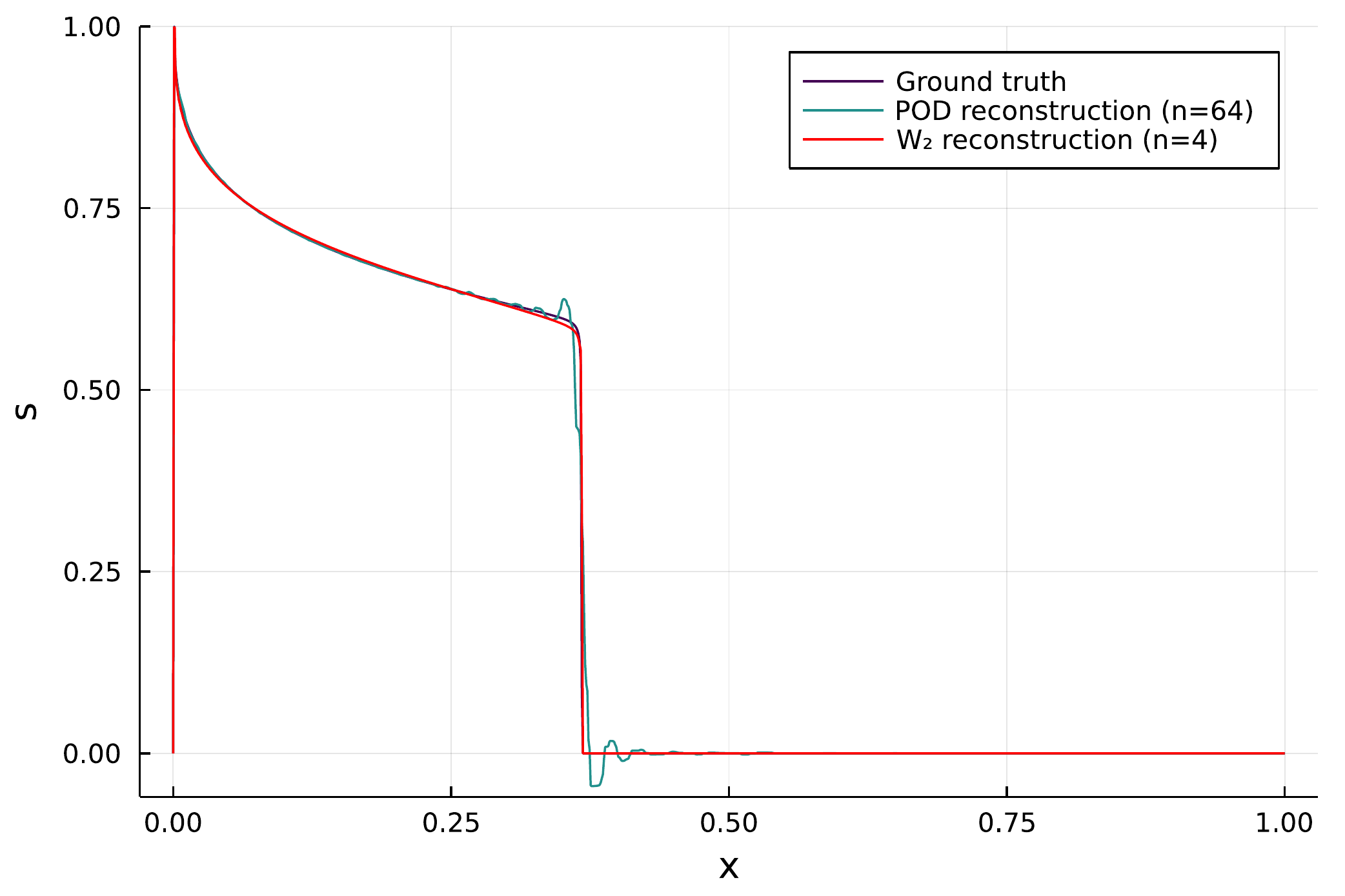}
        \caption{Example 1: Snapshot and its reconstruction using the gBar algorithm and a POD basis.}
        \label{fig:ex1_recon}
    \end{figure}

    Table \ref{tab:ex1_n_of_epsilon} shows that decent reconstructions can already be achieved with very few atoms. This highlights the strength of the method and its applicability to the example at hand. As expected, the POD reconstruction is ill-suited: to obtain an $L^1$-error of one percent, one needs over 230 modes, which is almost one third of the size of the training set - barely any reduction is possible. This is mainly due to oscillations occurring in the POD reconstruction around the jump discontinuity as shown in Figure \ref{fig:ex1_recon}.
    
    \begin{table}[h!]
    \begin{center}
        \begin{tabular}{ |c|c|c|c|c| } 
        \hline 
        $\varepsilon$ & 0.1 & 0.05 & 0.01 & 0.005 \\ 
        \hline
        $n_{gBar}(\varepsilon)$ & 2 & 3 & 16 & 29 \\ 
        $n_{POD}(\varepsilon)$ & 41 & 77 & 231 & $\gg$ 250 \\ 
        \hline 
        \end{tabular} \\ \medskip
        \end{center}
        \caption{Atoms and POD modes needed to reach a given tolerance $\varepsilon$, in $L^1$ norm, for the snapshots used in example 1.}
        \label{tab:ex1_n_of_epsilon}
        \end{table}
    
\subsection{Example 2: Porous medium made of two rock types with varying permeability and position of the heterogeneity}

\begin{figure}[htbp]
        \centering
        \includegraphics[width=0.5\textwidth]{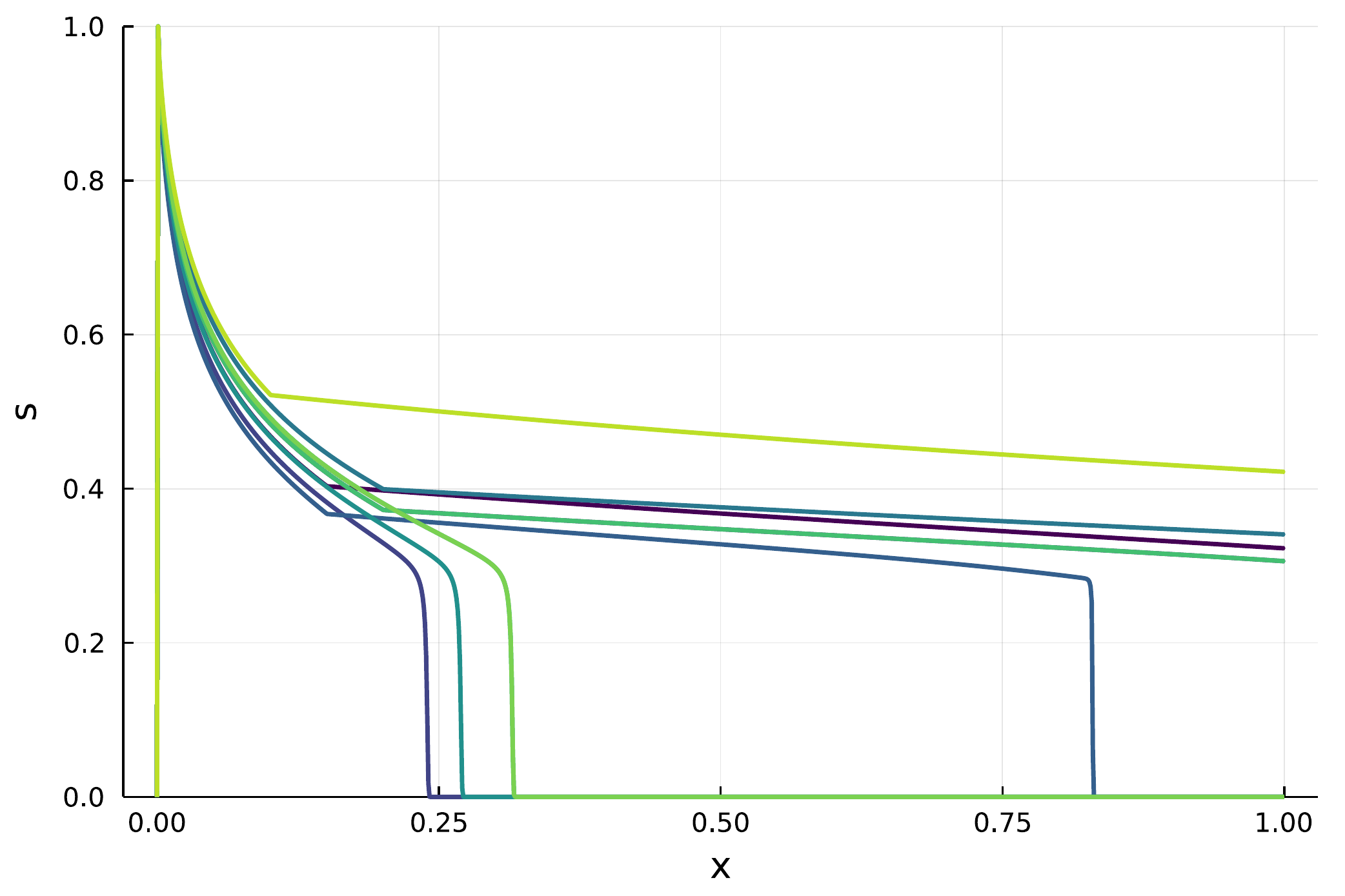}
        \caption{Randomly selected snapshots at final time from example 2.}
        \label{fig:ex2_snaps}
    \end{figure}

    For the second example, a heterogeneous medium is tested, composed of two different rock types: a highly permeable one (HP) characterized by $\phi_{HP}=0.1$ and $k_{HP}=10^{-13} \, \mathrm{m^2}$, and a low permeable one (LP) with fixed porosity $\phi_{LP}=0.01$ but uniformly varying permeability $k_{LP} \in \{ 5,6,7,8,9 \} \times 10^{-14} \, \mathrm{m^2}$. In this case, the viscosities are kept constant at $\mu_{w}=0.003 \, \mathrm{Pa\,s}$ and $\mu_{nw}=0.03 \, \mathrm{Pa\,s}$ (giving a ratio of $\mu=10$), and so is the relative permeability exponent, $\beta=2$.\\
    The location of the interface between the two rock types is denoted by $\gamma$. For $x<\gamma$, the rock type is HP, and LP otherwise. $\gamma$, which is sampled at irregular intervals, with a finer sampling around $x = 0$: $\gamma \in \{ 0.0 ,0.05,0.1,0.15,0.2,0.4,0.6 \}$. The evolution of the water saturation is calculated over 10 years and snapshots are regularly taken at $t \in \{0.5, 1.0, \dots, 10.0 \}$ years, leading to $z = (t,k_{LP},\gamma)$ and $| \cZ_{\rm train} | = 700$. \\ From Figure \ref{fig:ex2_snaps}, where some snapshots at the final time are shown, we see that the solution is quite different to the homogeneous case: when the saturation front reaches the permeability interface, the velocity of propagation changes, leading to a kink in the saturation profile. In some cases, the front does not reach the jump and the evolution of the saturation still behaves as in an homogeneous medium. \\
    
\begin{figure}[htbp]
    \centering
    \begin{subfigure}{.47\textwidth}
        \centering
        \includegraphics[width=0.92\textwidth]{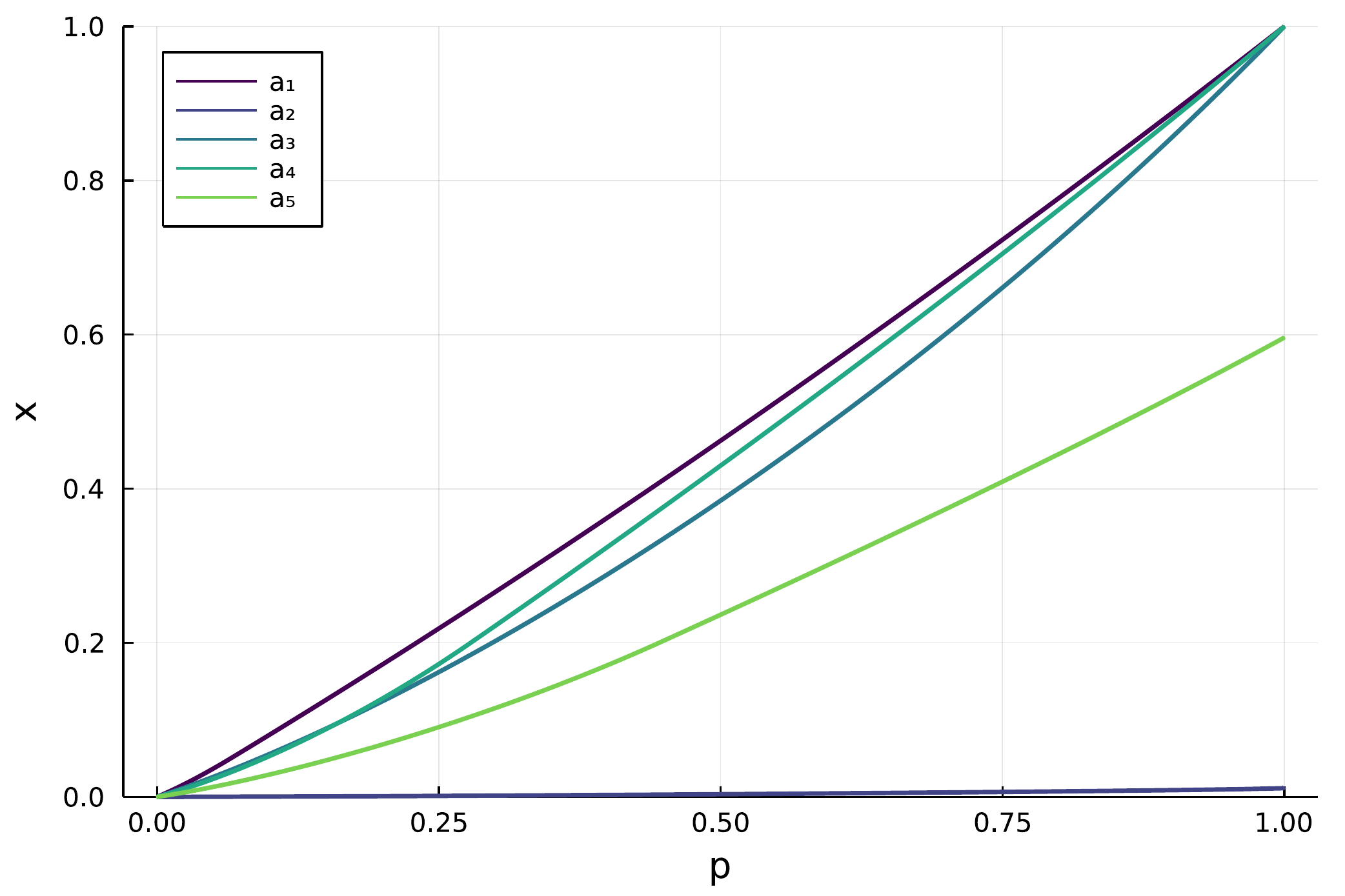}
        \caption{Dictionary atoms...}
    \end{subfigure}%
    \hspace{0.04\textwidth}
    \begin{subfigure}{.47\textwidth}
    \centering
        \includegraphics[width=0.92\textwidth]{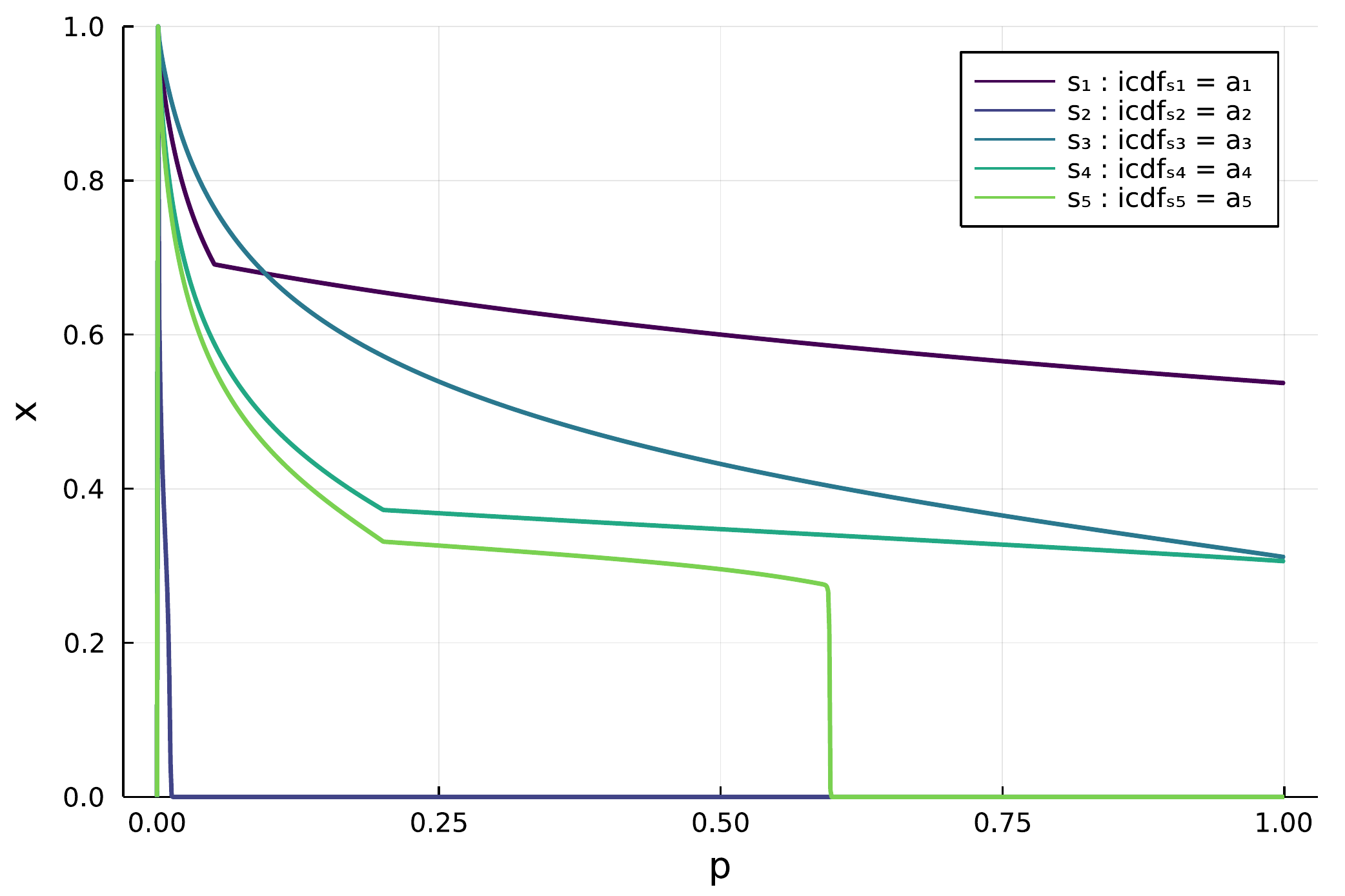}
        \caption{...and corresponding densities.}
    \end{subfigure}
    \caption{First five atoms chosen by the greedy algorithm for example 2. }
    \label{fig:ex2_atoms}
    \end{figure}

\begin{figure}[htbp]
    \centering
    \begin{subfigure}{.47\textwidth}
    \centering
    \includegraphics[width=0.92\textwidth]{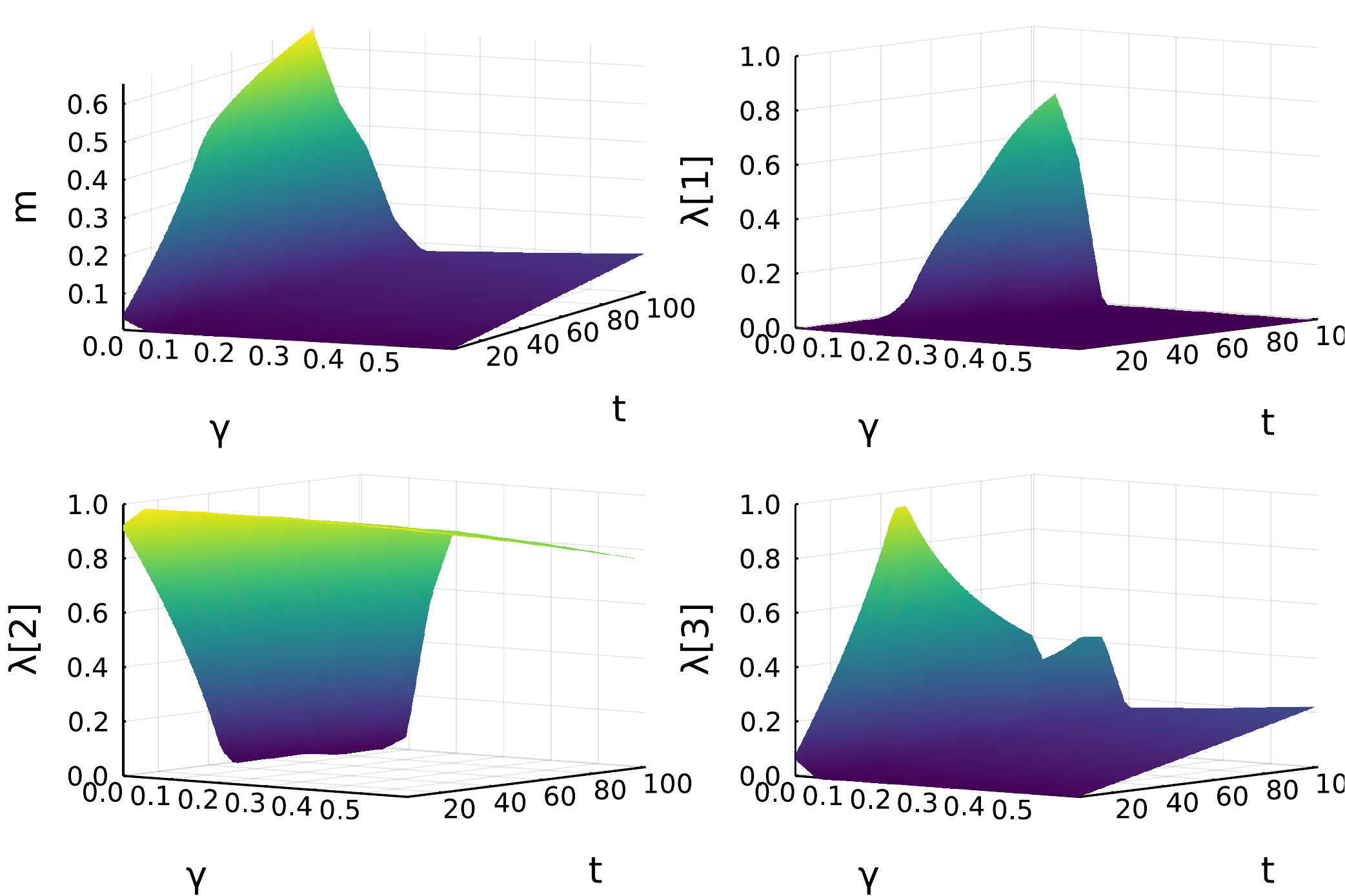}
    \caption{$(\gamma,t) \mapsto \lambda_{1 \leq i \leq 3}(t, \gamma, k_{LP} = 5 \times 10^{-14})$ and $(\gamma,t) \mapsto m(t, \gamma, k_{LP} = 5 \times 10^{-14})$.}
    \label{fig:ex2_interpolates1}
    \end{subfigure}%
    \hspace{0.04\textwidth}
    \begin{subfigure}{.47\textwidth}
    \centering
    \includegraphics[width=0.92\textwidth]{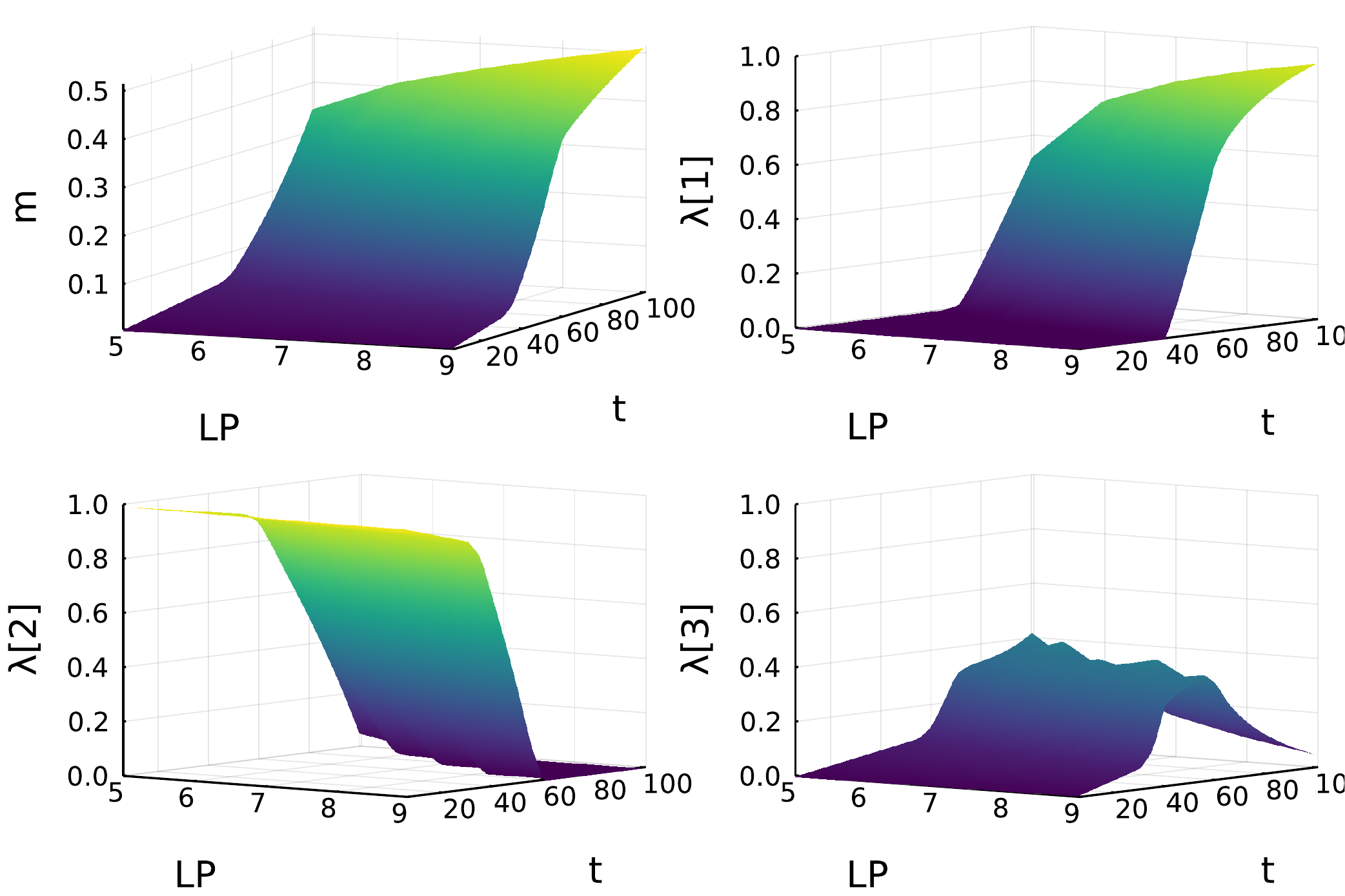}
    \caption{$(k_{LP},t) \mapsto \lambda_{1 \leq i \leq 3}(t, \gamma = 0.1, k_{LP})$ and $(k_{LP},t) \mapsto m(t, \gamma = 0.1, k_{LP})$.}
    \label{fig:ex2_interpolates2}
    \end{subfigure}
    \caption{Interpolations from parameter space to barycentric weights and mass in example 2.}
    \label{fig:ex2_interpolates}
    \end{figure}

    The first atoms selected by the greedy algorithm are shown in Figure \ref{fig:ex2_atoms} and cross-sections of the interpolating functions for $\lambda_{1 \leq i \leq 3}(t,\gamma,k_{LP})$ and $m(t,\gamma,k_{LP})$ in Figure \ref{fig:ex2_interpolates}.
    
\begin{figure}[htbp]
        \centering
        \includegraphics[width=0.5\textwidth]{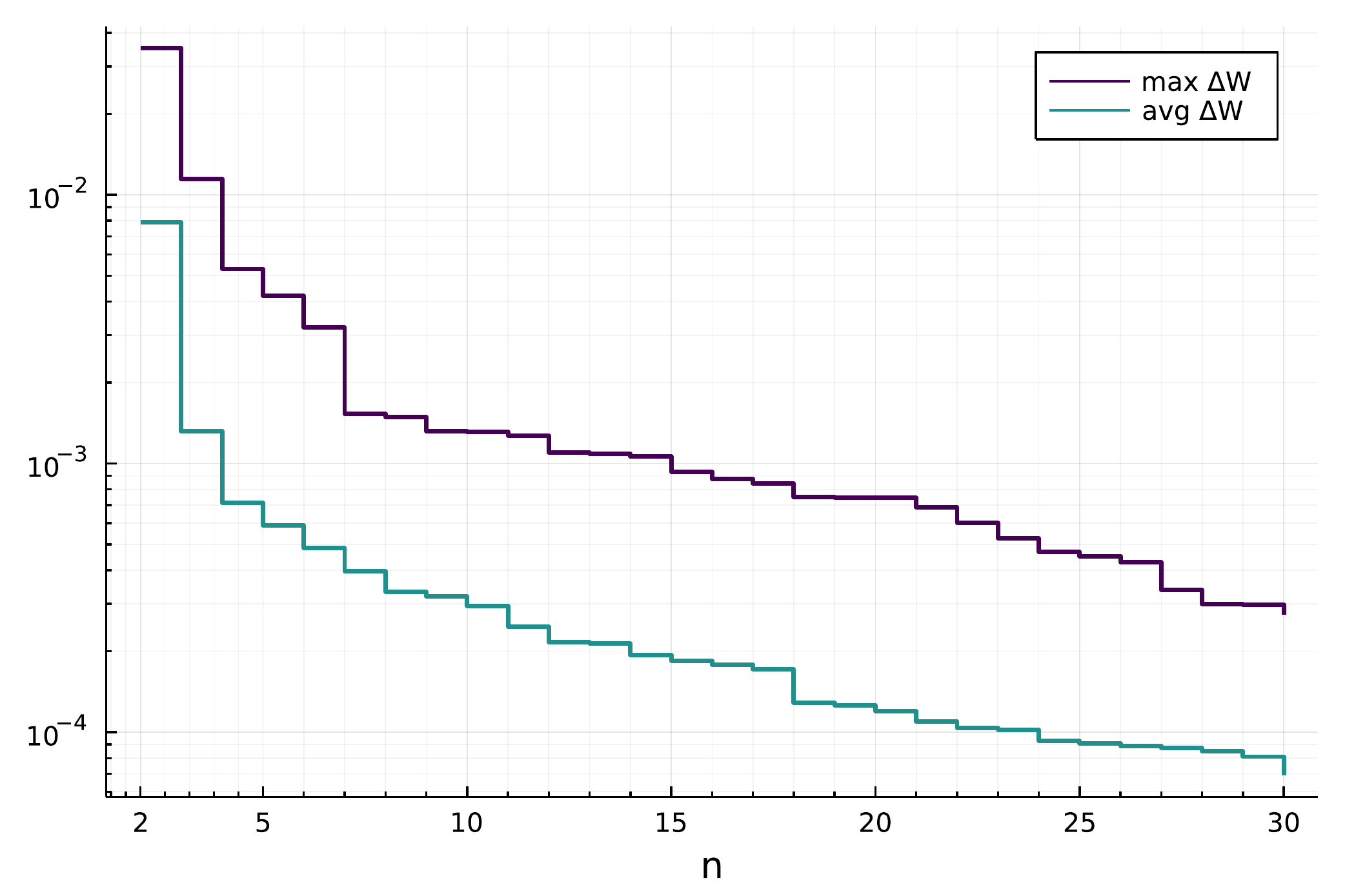}
        \caption{Decay of the $W_2$ error through the iterations of the greedy algorithm for example 2.}
        \label{fig:ex2_deltaW}
    \end{figure}

    Letting the greedy algorithm add atoms allows the $W_2$-error to decrease in a similar way with respect to the homogeneous case (Figure \ref{fig:ex2_deltaW}). Note, however, that both maximum and average errors in $W_2$ are about one order of magnitude above the values from Example 1. This also reflects in the errors as measured in the $L_1$ norm, shown in Table \ref{tab:ex2_n_of_epsilon}. It was not possible to obtain an accuracy of 0.5\% for this test-case before the barycentric algorithm became unstable and the iterative solver failed to converge to optimal barycentric weights.
    
    \begin{table}[htbp]
    \begin{center}
        \begin{tabular}{ |c|c|c|c|c| } 
        \hline 
        $\varepsilon$ & 0.1 & 0.05 & 0.01 & 0.005 \\ 
        \hline
        $n_{gBar}(\varepsilon)$ & 3 & 3 & 14 & - \\ 
        $n_{POD}(\varepsilon)$ & 20 & 35 & 118 & 173 \\ 
        \hline 
        \end{tabular} \\ \medskip
        \end{center}
        \caption{Atoms and POD modes needed to reach a given tolerance $\varepsilon$, in $L^1$ norm, for the snapshots used in example 2.}
        \label{tab:ex2_n_of_epsilon}
        \end{table}
        
\begin{figure}[htbp]
        \centering
        \includegraphics[width=0.5\textwidth]{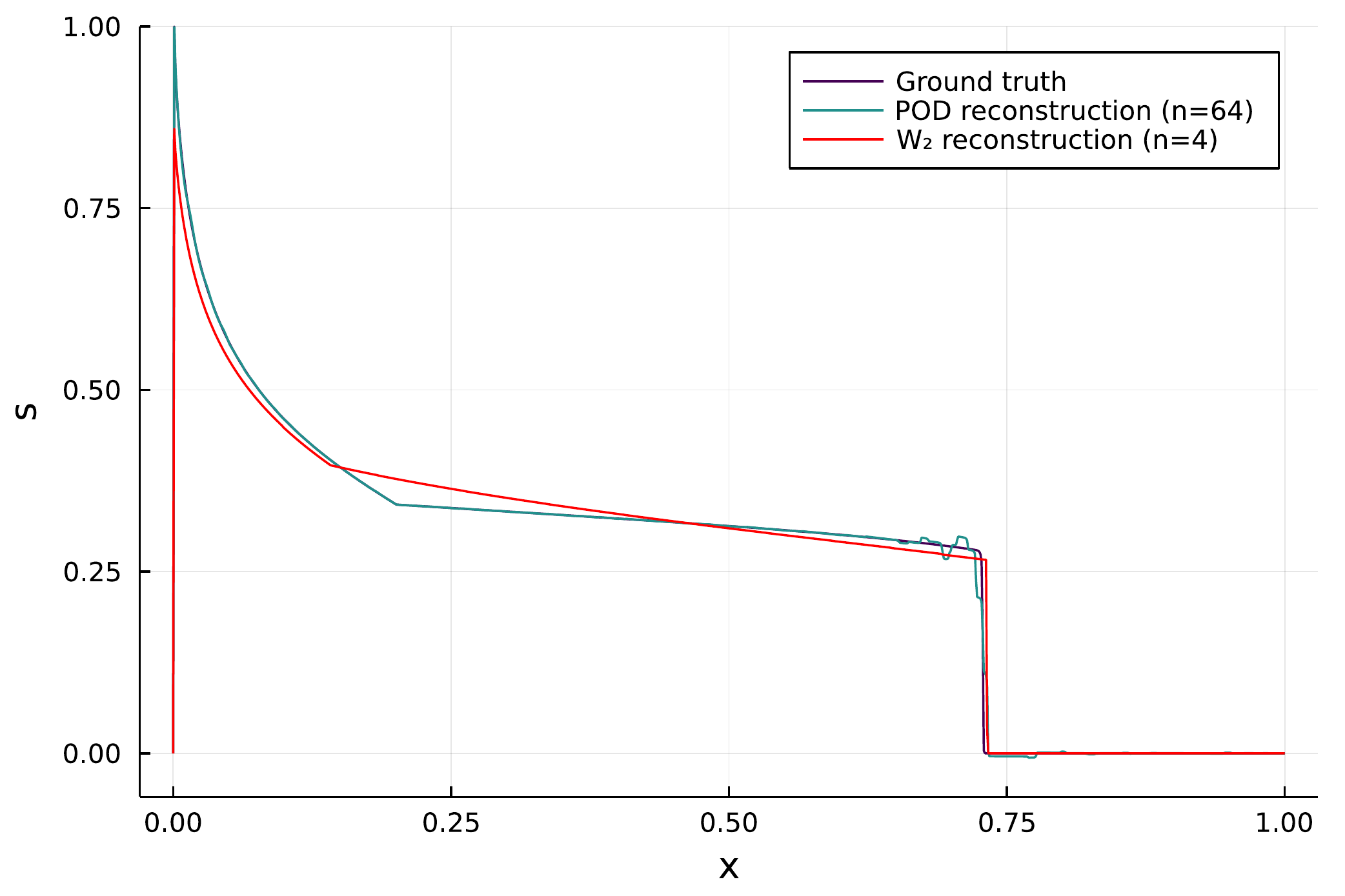}
        \caption{Example 2: Snapshot and its reconstruction using the gBar algorithm and a POD basis. Note how the barycentric interpolation does not capture the kink located at the HP/LP transition well.}
        \label{fig:ex2_recon}
    \end{figure}

    The reconstruction of a snapshot, plotted in Figure \ref{fig:ex2_recon}, shows that the proposed method does not accurately capture the kink of the solution at the interface of the HP/LP regions at $x \approx 0.2$. On the other hand, the POD reconstruction shows an unstability at the solution front. Several elements of the training set no longer show this jump discontinuity, however, and in these cases the POD method performs better, which reflects also in the second row of Table \ref{tab:ex2_n_of_epsilon}.
    
\subsection{Properties of the minimization problem}\label{sec:conditioning_convexity}
    
    \begin{figure}[htbp]
        \centering
        \includegraphics[width=0.5\textwidth]{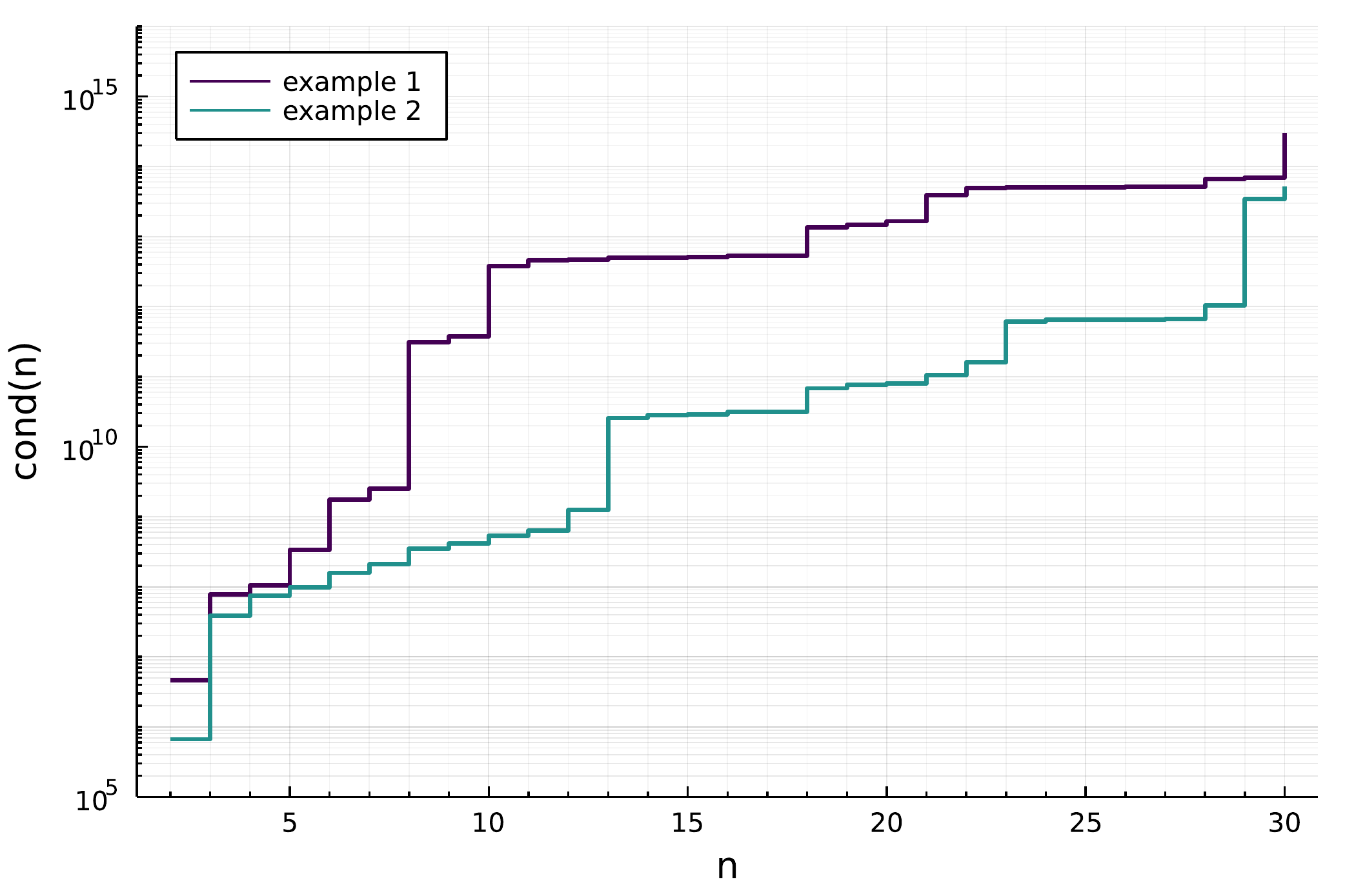}
        \caption{Conditioning number of the matrix $\mathtt{A^T A} $ as $n$ increases.}
        \label{fig:conditioning}
    \end{figure}
    
    In this section, we want to provide some insights into the properties of the optimization problem $\min_{\mathtt{\Lambda}} \Vert \mathtt{A \Lambda - f} \Vert$, i.e. the approximation of a discrete icdf $\mathtt f$ as a convex combination of a collection of atoms (the columns of $\mathtt A$). Figure \ref{fig:conditioning} shows the conditioning number of the matrix $\mathtt{A^T A}$ throughout the greedy algorithm iterations for examples one and two. In both cases, the conditioning number grows fast with increasing $n$, making it harder to solve for the optimal $\mathtt \Lambda$. Note that a bad conditioning of the matrix $\mathtt{A^T A}$ does not necessarily mean that the optimization problem has no solution. For example, for $s_1 = \delta_{x_1}$ and $s_2 = \delta_{x_2}$, we find $\icdf_1 \equiv x_1$ and $\icdf_2 \equiv x_2$ and $\mathtt{A^T A}$ is singular. Nevertheless, the optimization problem for the optimal barycenter to approximate a given $u$ by these atoms is well posed - it is solved by projecting the mean of $u$ onto $\{ (1-\lambda)x_1 + \lambda x_2 : 0 \leq \lambda \leq 1 \}$.
    
    \medskip

\begin{figure}[htbp]
    \centering
    \begin{subfigure}{.3\textwidth}
    \centering
    \includegraphics[width=\textwidth]{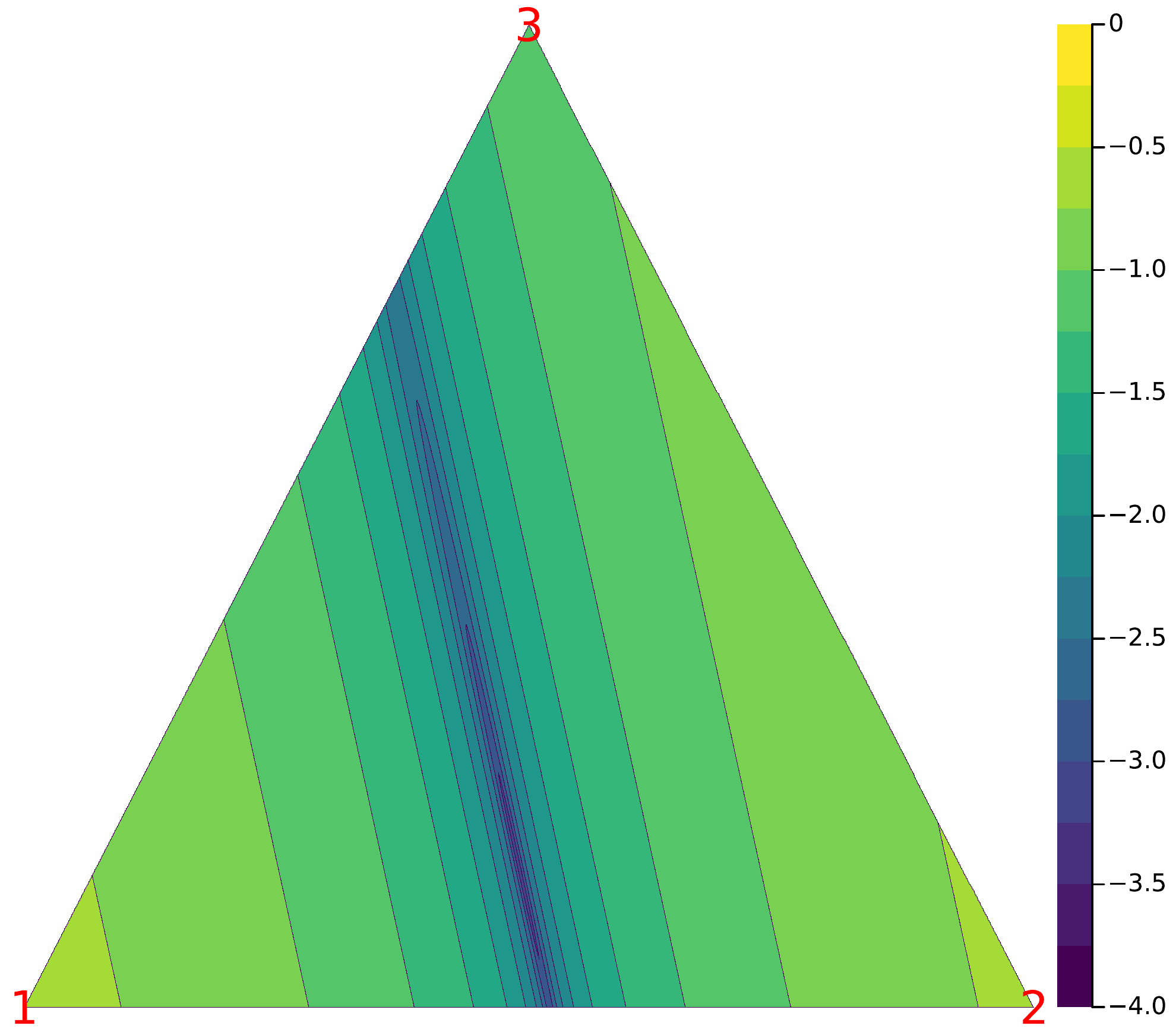}
    \caption{$\Lambda^{\rW opt}_3 \approx (0.44, 0.42, 0.14)$ and $\min \Delta W_2 \approx 7.8 \times 10^{-4}$}
    \label{fig:ex1_energy_landscape3}
    \end{subfigure}%
    \hspace{0.02\textwidth}
    \begin{subfigure}{.3\textwidth}
    \centering
    \includegraphics[width=\textwidth]{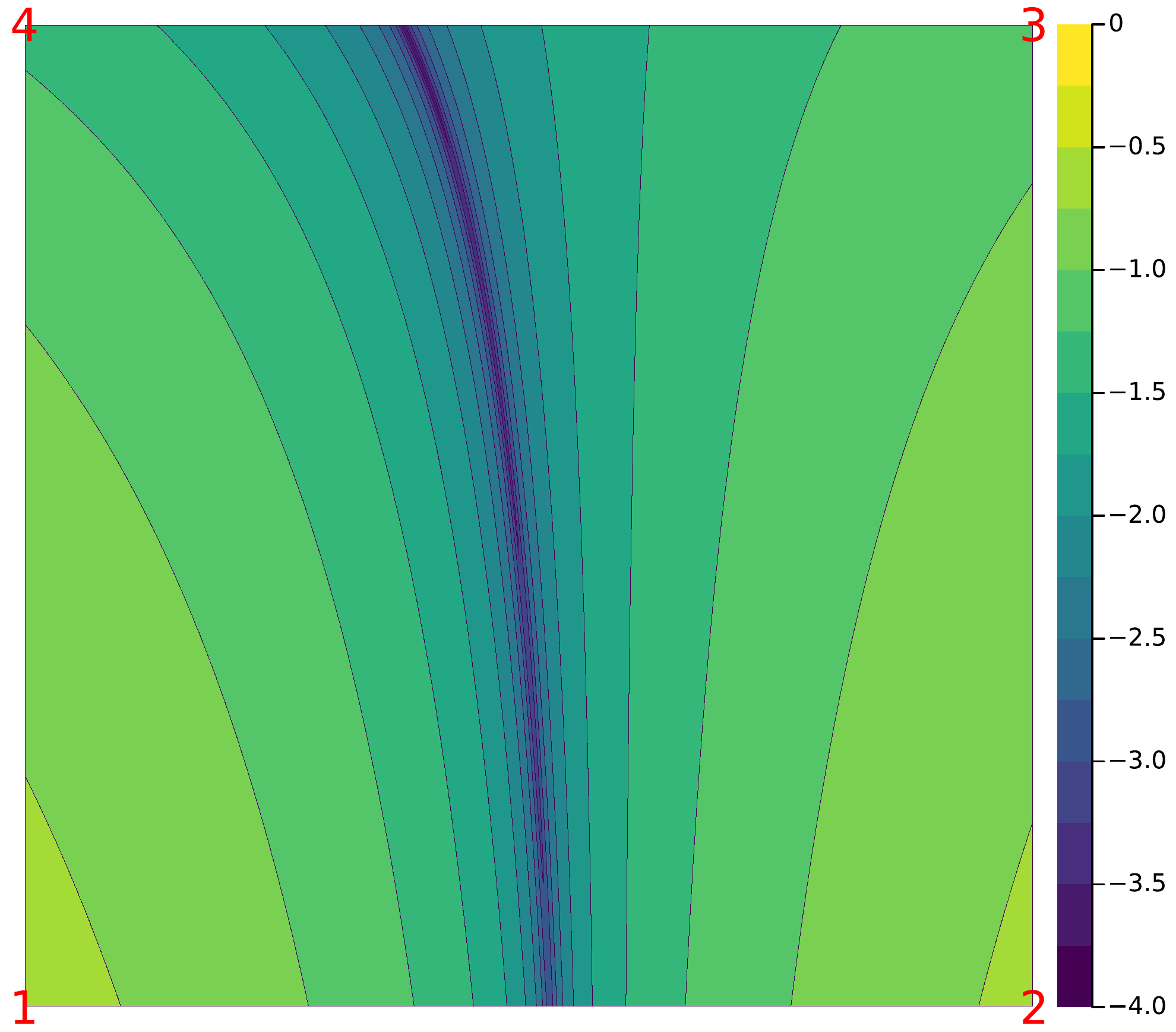}
    \caption{$\Lambda^{\rW opt}_4 \approx (0.01, 0.00, 0.38, 0.61)$ and $\min \Delta W_2 \approx 1.9 \times 10^{-4}$}
    \label{fig:ex1_energy_landscape4}
    \end{subfigure}%
    \hspace{0.02\textwidth}
    \begin{subfigure}{.3\textwidth}
    \centering
    \includegraphics[width=\textwidth]{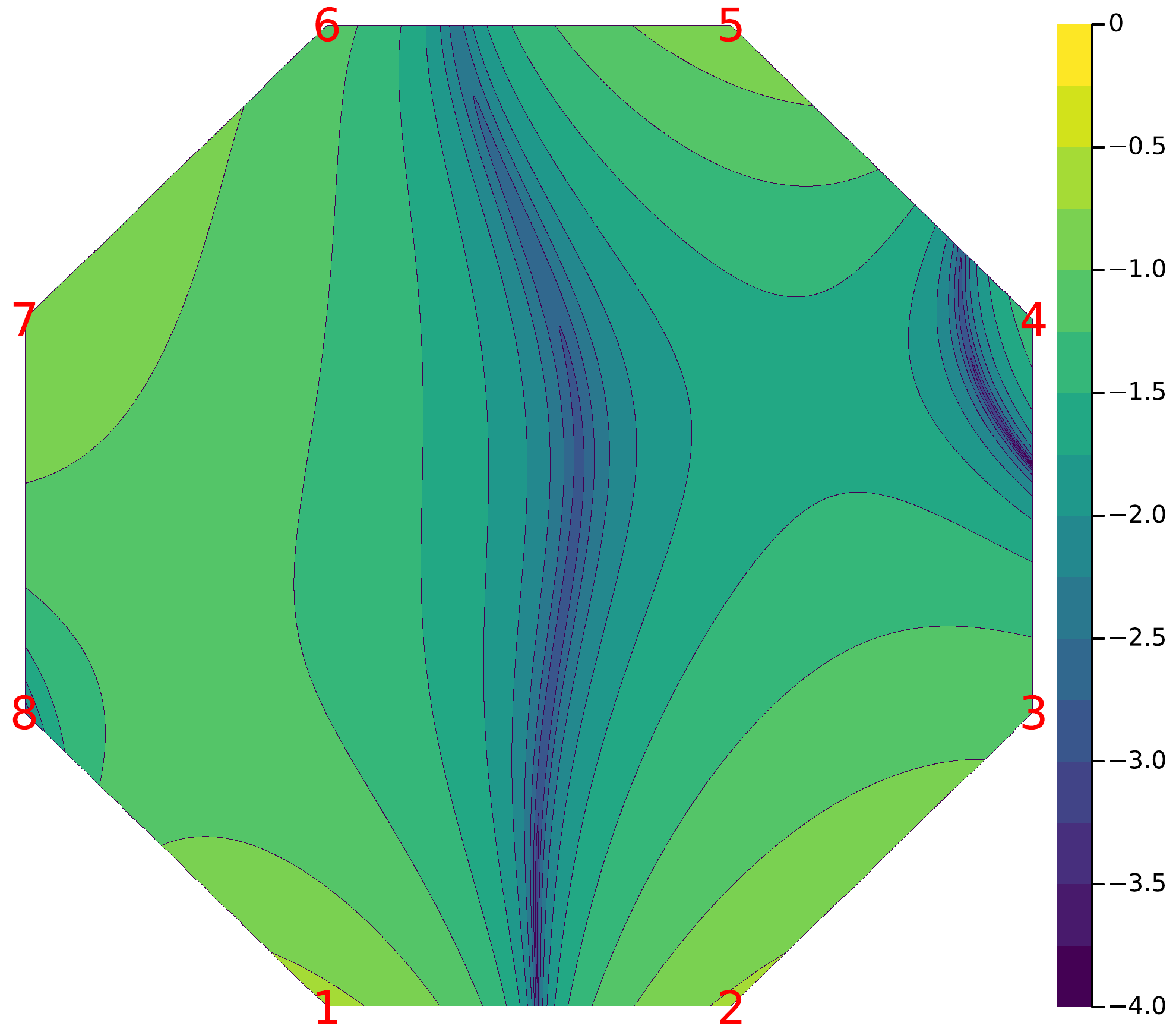}
    \caption{$\Lambda^{\rW opt}_8 \approx (0.0, 0.0, 0.0, 0.26, 0.42, 0.0, 0.31, 0.0)$ and $\min \Delta W_2 \approx 1.6 \times 10^{-4}$}
    \label{fig:ex1_energy_landscape8}
    \end{subfigure}%
    \\
    \begin{subfigure}{.3\textwidth}
    \centering
    \includegraphics[width=\textwidth]{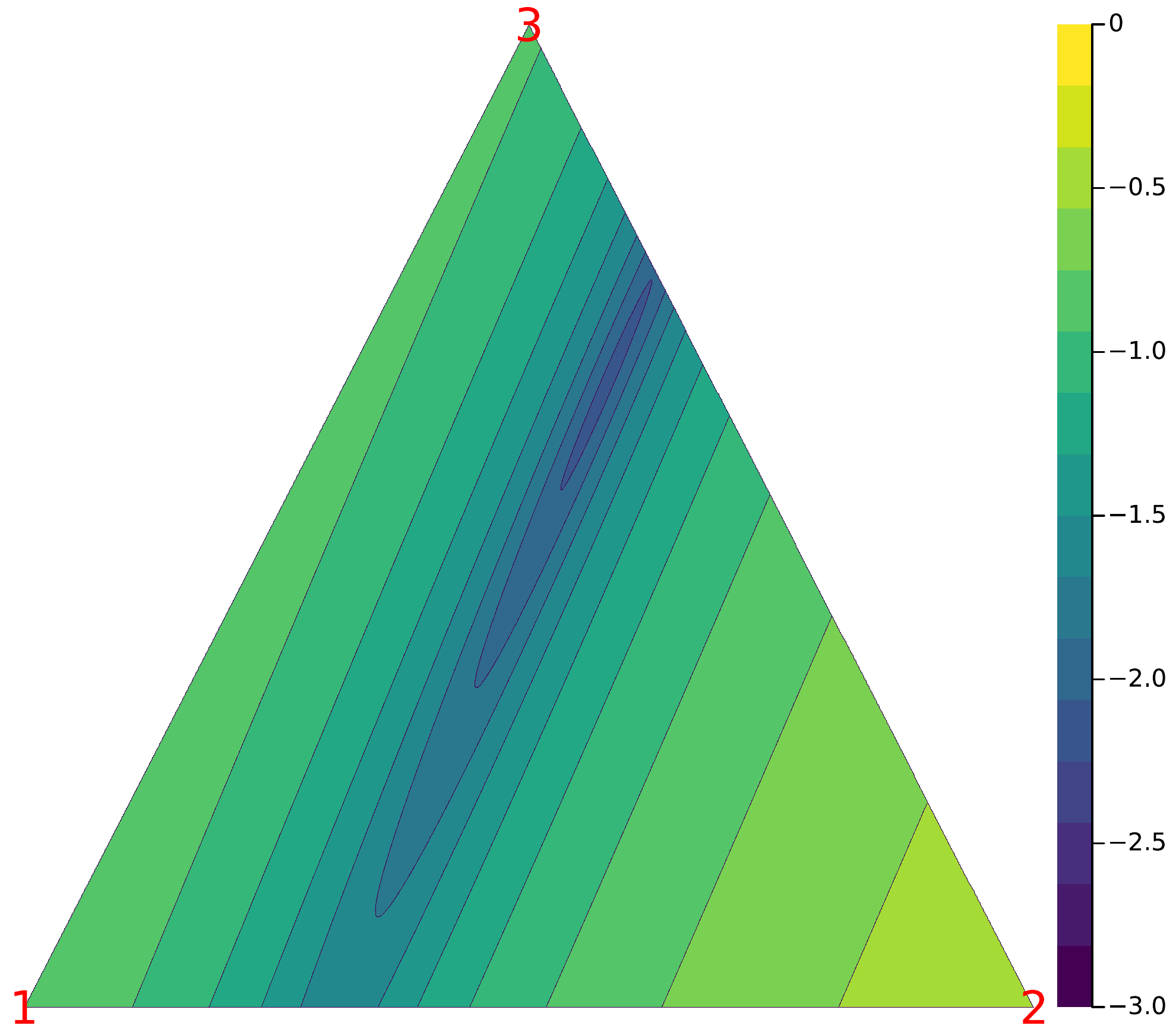}
    \caption{$\Lambda^{\rW opt}_3 \approx (0.11, 0.26, 0.63)$ and $\min \Delta W_2 \approx 7.8 \times 10^{-3}$}
    \label{fig:ex2_energy_landscape3}
    \end{subfigure}%
    \hspace{0.02\textwidth}
    \begin{subfigure}{.3\textwidth}
    \centering
    \includegraphics[width=\textwidth]{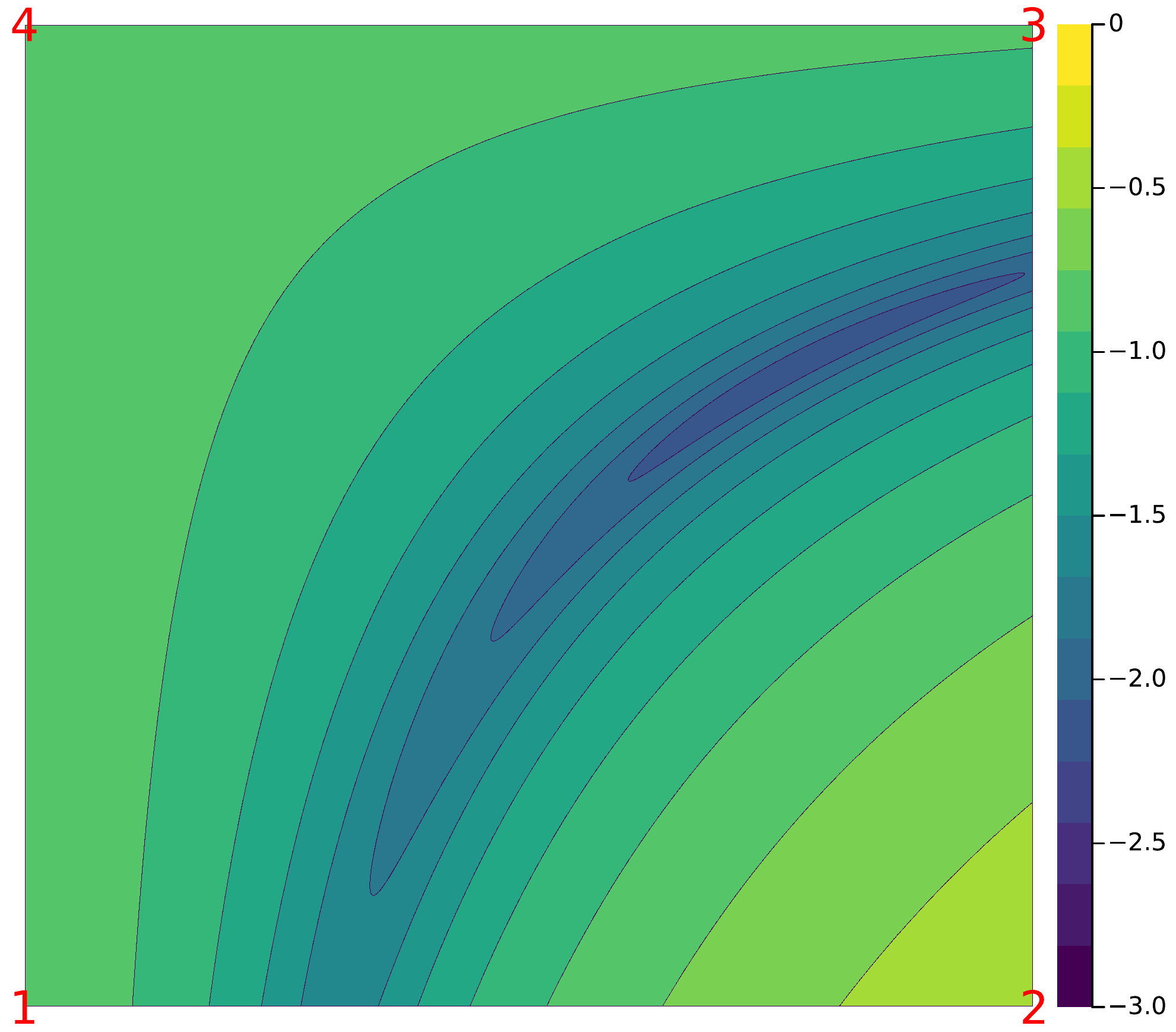}
    \caption{$\Lambda^{\rW opt}_4 \approx (0.0, 0.27, 0.33, 0.40)$ and $\min \Delta W_2 \approx 4.7 \times 10^{-3}$}
    \label{fig:ex2_energy_landscape4}
    \end{subfigure}%
    \hspace{0.02\textwidth}
    \begin{subfigure}{.3\textwidth}
    \centering
    \includegraphics[width=\textwidth]{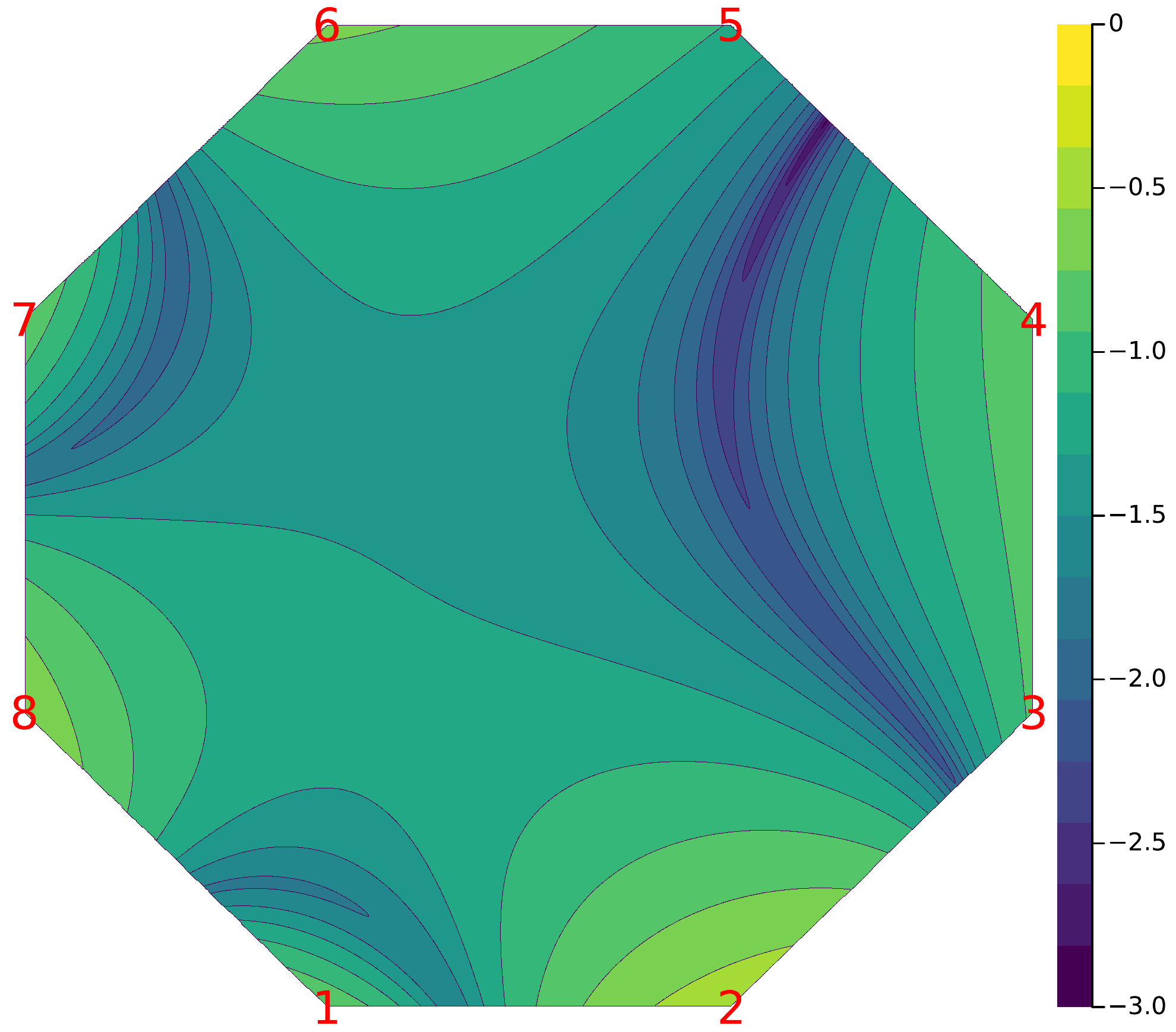}
    \caption{$\Lambda^{\rW opt}_8 \approx (0.0, 0.0, 0.03, 0.30, 0.67, 0.0, 0.0, 0.0)$ and $\min \Delta W_2 \approx 1.3 \times 10^{-3}$}
    \label{fig:ex1_energy_landscape8_2}
    \end{subfigure}%
    \caption{$\log_{10}$ of the $W_2$ distance between barycenter reconstructions and the snapshot from Figure \ref{fig:ex1_recon} (top row) and Figure \ref{fig:ex2_recon} (bottom row) for different $n$ as a function of $\lambda_1, \dots, \lambda_n$. The vertices of the polygons mark the atoms $a_1, \dots, a_n$.}
    \label{fig:energy_landscapes}
    \end{figure}

    Figure \ref{fig:energy_landscapes} shows the $W_2$ error between two target icdfs and the barycenter of atoms for different $n$ and $\Lambda$ in log-scale. We display elements of the $n-1$-dimensional simplex as points in a convex polygon. For $n=3$, this corresponds to the well-known barycentric coordinates within a triangle, depicted in Figures \ref{fig:ex1_energy_landscape3} and \ref{fig:ex2_energy_landscape3}. For $n>3$, we use generalized barycentric coordinates (in particular, so-called Wachspress coordinates) within polygons as described in \cite{floater_generalized_2015}. These coordinates associate to every point $x_{\rm poly} \in \{ n-\text{Polygon} \} \subset \mathbb R^2$ a vector of barycentric weights $\lambda_{\rm gbc} \in \Sigma_n \subset \mathbb R^n$. The $i$-th vertex of a polygon corresponds to barycentric weights with only one non-zero entry for the $i$-th atom, while points along the edge between the vertices corresponding to $a_i$ and $a_j$ represent $(1-\lambda)a_i + \lambda a_j$ where $0 < \lambda < 1$ along the edge. The center of each polygon corresponds to the point where $\lambda_{{\rm gbc}, i} = 1/n \; \forall i$. The mapping $x_{\rm poly} \mapsto \lambda_{\rm gbc}$ is not surjective, so the second and third column of Figure \ref{fig:energy_landscapes} show only a specific``cut" through the true energy landscape defined on $\Sigma_4$ and $\Sigma_8$. For example, from Figure \ref{fig:ex1_energy_landscape4}, we can see that approximations with $W_2$ error below $10^{-3}$ lie close to a line that enters the tetrahedron $\Sigma_4$ at the $a_1$-$a_2$ edge and exits at the $a_3$-$a_4$ edge. The minimum lies on the $a_1$-$a_3$-$a_4$ face, very close to the $a_3$-$a_4$ edge. In the case $n=8$, it naturally is quite hard to faithfully represent this high-dimensional setting with a two-dimensional plot. Recall that the objective is a quadratic function constrained to a convex set. Therefore, no disjoint minima are possible. 
    
    In all depicted cases, the contour lines of the energy are very eccentric, a consequence of the poor conditioning of the matrix $\mathtt{A^T A}$. In both examples and for large $n$, optimal barycentric weights lie on the faces of the simplices, with several zero entries. This sparsity is seen throughout the elements of $\cM_{\rm train}$ and has also been observed in \cite{bonneel_wasserstein_2016}. It is an interesting question for future work how one could choose and/or modify the atoms in order to avoid this and thereby moving the minimum of the objective into the interior of the simplex. 
    
    \medskip
    
\begin{figure}[htbp]
        \centering
        \includegraphics[width=0.5\textwidth]{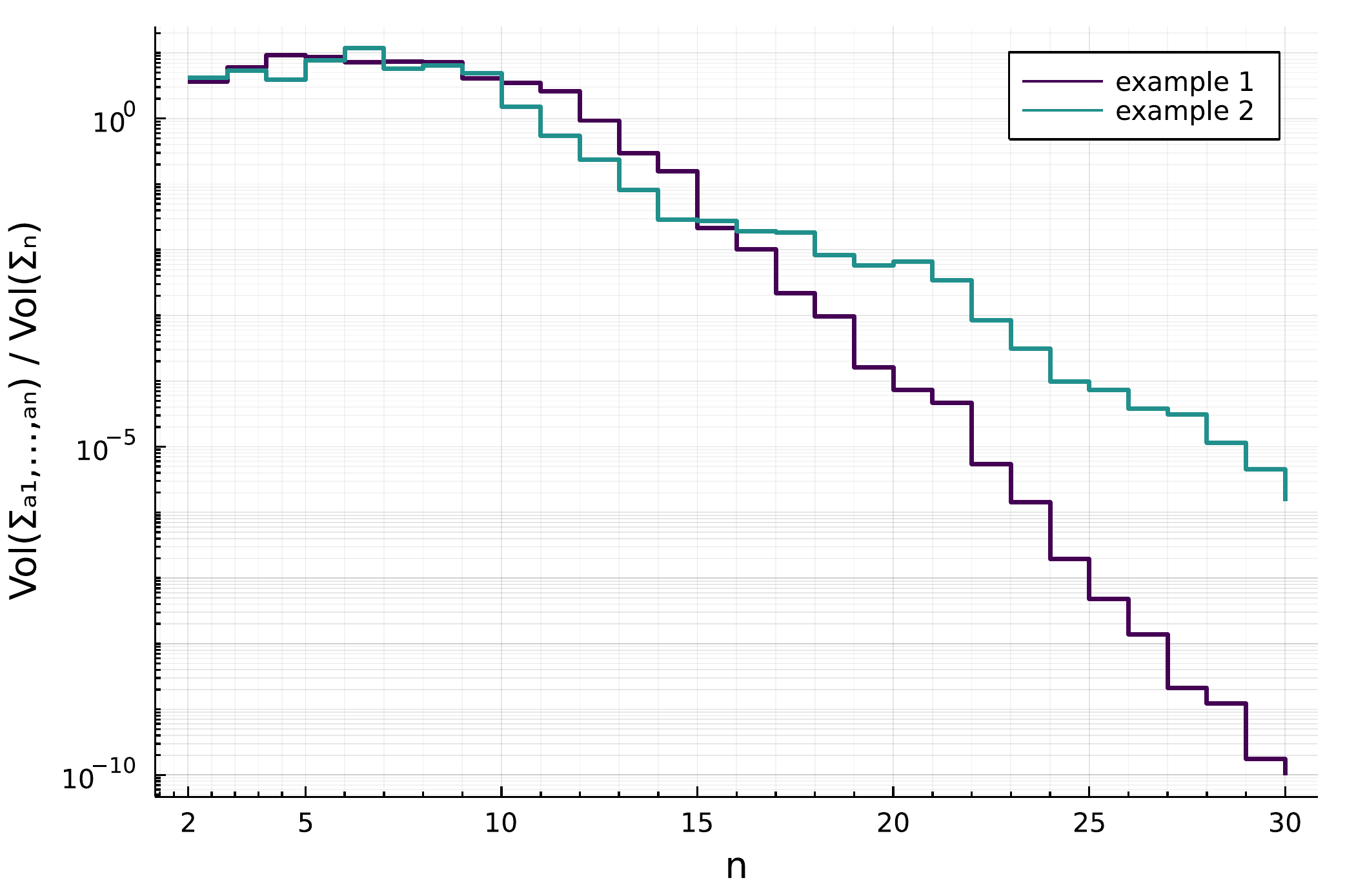}
        \caption{Normalized volume of the simplex with vertices $a_1, \dots, a_n$ as $n$ increases.}
        \label{fig:simplexvolume}
    \end{figure}

    The existence of a mapping from the Wasserstein space to $\mathbb R^N$ (namely, $\pdf_s \mapsto \icdf_s$) is a special case of one spatial dimension and does not hold in higher-dimensional case (\cite{peyre_computational_2020}, Proposition 8.2). To get an understanding if we are adding an atom to the collection in step $n+1$ that is already close to a convex combination of the atoms $\{ a_1, \dots, a_n \}$, we can exploit this. Using the distance matrix $\left (  \Vert a_i - a_j \Vert^2_{L_2} \right )_{1 \leq i,j \leq n}$, we calculate the volume of the simplex with vertices $\{a_i\}_{1 \leq i \leq n}$ as a function of $n$ using the Cayley–Menger determinant~\cite{doi:10.4169/amer.math.monthly.124.7.621}. We normalize this value with the volume of the unit $n$-simplex. The results are shown in Figure \ref{fig:simplexvolume}. 
    
    Adding another atom that is close to the simplex with vertices $a_1, \dots, a_n$ will lead to a small volume for the resulting simplex with vertices $a_1, \dots, a_{n+1}$. For $n=2$, this would be a very slim triangle. We observe that up until $n\approx 10$, in both examples, the simplices relative volume remains almost constant, which we interpret as every additional atom providing incremental information. For larger $n$, the relative volume decreases exponentially as the simplices become more degenerate, i.e. some of the $a_i$'s are almost affinely dependent.
    
\subsection{Conclusions from the numerical experiments}
    
    The numerical experiments exemplify that the presented method is able to represent the solution space of a complex non-linear equation using only very few atoms. In the first test-case considered, using only $n=3$ atoms allows to reconstruct every element in training data set with an accuracy of 95\%. The data depend on three parameters and the solutions vary considerably with those parameter values - note, for example, how slowly the front of the wetting phase moves at high values of $\mu$ and $\beta$. The map from the parameter space to the barycentric weights for only a small number of atoms has a simple structure as seen in Figure \ref{fig:ex1_interpolates1} and \ref{fig:ex1_interpolates2}. \\
    The second test-case features a collection of solutions that are both more challenging for the presented method and easier for the POD method to represent. In this case, several elements of $\mathcal{M}_{\rm train}$ are supported on the entire domain, which removes the biggest challenge for the POD representation, i.e. the jump discontinuity. At the same time, the barycentric interpolation does not preserve the characteristic kink located at the HP/LP interface, reducing the quality of the barycentric reconstruction. Ultimately, the barycentric method still performs well even with very few atoms, but cannot reach the same level of accuracy from test-case 1 before the dictionary becomes too ill-conditioned and the greedy algorithm terminates. \\
    It is important to note, when comparing the presented method to POD, that there is no straightforward online method for the latter. We are merely comparing the encoding and reconstruction properties of the two methods. The online method for the barycentric representation is purely interpolating and comes at a very small computational cost (that is: evaluation of the functions $\mathcal{I}_{1, \dots, n}$, calculation of the convex combination, inversion and differentiation of the icdf) that is linear in $N$, the dimension of the full-order problem.

\section{Future research directions and extensions}\label{sec:discussion}

    The computation of Wasserstein distances and barycenters through the inverse cumulative distribution functions is a special case limited to one-dimensional problems. There has been a lot of progress in the development of algorithms for arbitrary dimensions in recent years which can be used to extend the present work for these cases \cite{peyre_computational_2020}, \cite{sejourne_sinkhorn_2021}. We intend to address this in upcoming work.
    It can be worthwhile to allow for the inclusions of atoms into the dictionary that are not snapshots of the solution themselves. This approach is taken in \cite{schmitz_wasserstein_2018}, where the atoms are initialized as uniform distributions and then optimized alongside the barycentric weights, i.e. one solves
    
    \begin{equation}
    \min_{\rU_n \in \mathcal{P}_2(\Omega)^n} \min_{\Lambda_n \in \Sigma_n} \sum_{u \in \mathcal{U}_{\rm train}} W^2_2(u, \bary(\rU_n, \Lambda_n)).
    \end{equation}
    
    We would further like to explore ideas on how to improve the conditioning of the matrix arising from the atoms. Instead of adding new atoms to the dictionary as they come, one could look for a way of adding a modified atom which is - in a sense that has to be made precise - fully outside the current range of the dictionary. In classical linear reduced basis methods, this is ensured by projecting any new element to the orthogonal complement of the present basis but such an operation is not appropriate for the current framework. Greedy algorithms developed for the cone of positive functions such as the Enlarged Nonnegative Greedy algorithm introduced in \cite{BBMM2021} could be an interesting starting point for our purposes.

\section*{Funding and Acknowledgments}
    This work was initiated and substantially realized during the six weeks of the CEMRACS 2021 event held at CIRM, Luminy. The authors' stay at CIRM has been funded by IFPEN, ANR ADAPT (ANR-18-CE46-0001), ANR project
    COMODO (ANR-19-CE46-0002), and the Emergences Project grant ``Models and Measures'' from the Paris City Council. The authors would like to thank S\'ebastien Boyaval for fruitful discussions.

\bibliographystyle{siamplain}
\bibliography{references}

\end{document}